\documentclass[11pt]{article}   
\usepackage{amssymb,amscd,latexsym}   
\usepackage{amsmath}
\usepackage{amsthm}
\newcommand{\rar}{\rightarrow}
\newcommand{\lar}{\longrightarrow}
\newcommand{\llar}{-\kern-5pt-\kern-5pt\longrightarrow}

\newtheorem{Theorem}{Theorem}[section]
\newtheorem{Lemma}[Theorem]{Lemma}
\newtheorem{Corollary}[Theorem]{Corollary}
\newtheorem{Proposition}[Theorem]{Proposition}
\newtheorem{Remark}[Theorem]{Remark}
\newtheorem{Example}[Theorem]{Example}
\newtheorem{Conjecture}[Theorem]{Conjecture}
\newtheorem{Definition}[Theorem]{Definition}
\newtheorem{Question}[Theorem]{Question}
\def\sqr#1#2{{\vcenter{\hrule height.#2pt
        \hbox{\vrule width.#2pt height#1pt \kern#1pt
            \vrule width.#2pt}
        \hrule height.#2pt}}}
\def\phi{\varphi}
\def\demo{\noindent{\bf Proof. }}
\def\square{\mathchoice\sqr64\sqr64\sqr{4}3\sqr{3}3}
\def\qed{\hspace*{\fill} $\square$}

\def\xx{{\bf x}}

\def\tt{{\bf t}}

\def\ff{{\bf f}}

\def\ff{{\bf f}}

\def\hh{{\bf h}}

\def\fm{{\mathfrak m}}

\def\hht{{\rm ht}\,}

\def\restr{{\kern-1pt\restriction\kern-1pt}}

\def\pp{{\mathbb P}}

\begin{document}
\begin{center}
{\Large{\bf\sc Homaloidal nets and ideals of fat points II: subhomaloidal nets}}
\footnotetext{Mathematics Subject Classification 2010
 (MSC2010). Primary  13D02, 13H10, 13H15, 14E05, 14E07, 14M05;
Secondary 13A02,  13C14,  14C20.}

\vspace{0.3in}

\hspace{-15pt}{\large\sc Zaqueu Ramos}\footnote{Partially
supported by a CNPq post-doc fellowship (151229/2014-7).}
\quad\quad\quad
 {\large\sc Aron  Simis}\footnote{Partially
supported by a CNPq grant (302298/2014-2) and by a CAPES-PVNS Fellowship (5742201241/2016).}

\end{center}


\bigskip

\begin{abstract}

This paper is a natural sequel to \cite{fat1} in that it tackles problems of the same nature. Here one aims at the ideal theoretic and homological properties of a class of ideals of general plane fat points whose second symbolic powers hold virtual multiplicities of proper homaloidal types.
For this purpose one carries a detailed examination of their linear systems at the initial degree, a good deal of the results depending on the method of applying the classical arithmetic quadratic transformations of Hudson--Nagata.
A subsidiary guide to understand these ideals through their initial linear systems has been supplied by questions of birationality with source $\pp^2$ and target higher dimensional spaces. This leads, in particular, to the retrieval of birational maps studied  by Geramita--Gimigliano--Pitteloud, including a few of the celebrated Bordiga--White parameterizations.

\end{abstract}

\section*{Introduction}

The title of this paper tells that it is a second part to \cite{fat1}, yet its contents are not  a logical sequel to the latter.
It would be fare to say that, while in \cite{fat1} the emphasis was upon the ideal theoretic questions of a proper homaloidal type, the present work stresses the behavior of other fat point ideals related to the latter.

Ideals of plane fat points have been thoroughly looked at in the last thirty years or so, both from the viewpoint of their Hilbert function and of the corresponding divisors on the blowing-up of $\pp^2$ at the given points.
While the appropriate geometric tool often works wonder on the blown-up territory, it is the algebraic side of the Hilbert function method that facilitates the use of plane Cremona transformations as a helpful tool.

To stay concrete, consider an ideal of fat points 
\begin{equation*}
J=I(\mathbf{p};{\boldsymbol\mu}):=\bigcap_{i=1}^r P_i^{\mu_i},
\end{equation*} 
where $P_i\subset R:=k[x,y,z]$ denotes the homogeneous prime ideal of the point $p_i\in \pp^2$, with $k$ denoting the ground field.
We assume throughout that the base points $p_i$ are general, a matter of no triviality.
The main fat ideal considered here comes out of requiring that its virtual multiplicities $\boldsymbol\mu$ satisfy a pair of equations in the same spirit of the equations of conditions satisfied by the virtual multiplicities of a homaloidal type.
In fact, the virtual multiplicities of the second symbolic power of such an ideal fulfill the equations of condition -- for this reason, and  in the lack of a better designation, this ideal is henceforth called  of {\em sub-homaloidal type}.

The main goal is to understand not just the Hilbert function of $J$ but also the nature and properties of its homogeneous parts (here referred by abuse as linear systems).
Following one of the main questions about fat points, we ask  when some of these linear systems attain their expected dimension.
There is a famous general conjecture in this connection, the so-called {\em Segre conjecture} -- nowadays known as the (alphabetically) {\em Gimigliano--Harbourne--Hirschowitz--Segre conjecture} after it has been realized that the original conjecture is equivalent to others considered by the first three authors.
Although this is usually phrased in terms of divisors and linear systems as global sections of rank one sheaves, we believe it becomes also transparent in the following form (see, e.g., \cite[Conjecture 1.2]{CiMi})

{\sc Conjecture.} If the general form of the $k$-vector space $J_d$ is reduced then $J_d$ has the expected dimension.

Simple as it sounds, this criterion is nevertheless hard to meet in theory.
It resembles the simple and equally hard  Enriques criterion for a linear system, satisfying the equations of condition, to define a Cremona map  in terms of the irreducibility of its general curves (see, e.g., \cite[Theorem 5.1.1]{alberich}).
 
One of the side results of our approach says that, for certain classes of sub-homaloidal fat ideals $J$, the linear system $J_d$ satisfy both the hypothesis and the conclusion of the above conjecture, when $d$ denotes the initial degree of $J$.
This knowledge does not logically produce a proof of the conjecture for such classes, but it reveals that the conjecture is accommodated in this setup. 

At the other end, the ideal of fat points $J$ treated in this paper, as well as in \cite{fat1}, lives on the far side with respect to the hypotheses of \cite{H0} 
 (see \cite{CiMi} for additional clarification) which require the traditional three highest virtual multiplicities to be small (a {\em standard} system) or else that the smallest virtual multiplicities be large (\cite[Remark 1.8]{CiMi}) with respect to the degree $d$ of the linear system $J_d$ in consideration.
 And indeed, the ideal $J$ studied here as well as in \cite{fat1} are typically, or closely related to, fat ideals defined by the base points and virtual multiplicities of plane Cremona maps.
 This implies that the corresponding linear system $J_d$ on the initial degree $d$ of $J$ satisfies the so-called Noether's inequality
 $$\mu_1+\mu_2+\mu_3 >d,$$
 where $\mu_1\geq \mu_2\geq \mu_3$ are the three highest virtual multiplicities of $J$.
 Linear systems with this condition have been dubbed {\em Cremona like} in \cite{CiMi}, while a standard linear system satisfies the reverse inequality.
We dare to say that most treatable classes of plane Cremona maps will have pretty small lowest virtual multiplicities even if the highest multiplicity is large (see \cite[Section 3]{fat1}).

 A marked feature of the ideals of fat points dealt with here or others  obtained from the latter by a quadratic transformation is that their initial degree $d$ and $\dim_k(J_d)$ are closely related.
 An additional important feature is the tendency to be minimally generated in at most two subsequent degrees -- a consequence of its Hilbert function being minimal. This is strictly so in the Cremona case (\cite[Theorem 2.6]{fat1}) and in the classes of sub-homaloidal ones mentioned above, as well as their quadratic transforms.
 Thus, the corresponding minimal free resolutions stay within the known conjectures found in the current literature.
 
 Finally, while a great many papers on the subject of linear systems of ideals of fat points ask when the linear system is empty, the present focus is on the linear system at the initial degree, hence taking non-emptiness for granted. The emphasis is accordingly moved to  the algebraic properties of the ideal, such as inquiring as to when it is generated in one single degree and, in this case, as to when it defines a birational map of $\pp^2$  onto the image and what properties the latter has.
 
 \smallskip
 
Having placed the object of the paper within the known literature, for the reader's convenience we now  confront it with the results of \cite{fat1}.
In the first three sections of \cite{fat1} we discussed the relation between the numerical and homological invariants of the {\em base ideal} of a plane Cremona map (i.e., the ideal of $R$ generated by its defining net of forms) and the ideal of its base points and virtual multiplicities.
To keep the discussion at the level of $\pp^2$ (without introducing blowing-up devices), we assumed throughout that the base points are proper, i.e., no infinitely base points around.
This is not a terrible shortcoming because by and large we need the points to be general to make sure they are base points of a Cremona map.
In the fourth (last) section of \cite{fat1} we reversed somehow the focus by starting up with  an ideal $J\subset R$ of fat points with multiplicities satisfying a couple of arithmetical conditions implying that the second symbolic power $J^{(2)}$ is the ideal of fat points coming from a {\em homaloidal type}, i.e., a vector of virtual multiplicities fulfilling the classical equations of conditions in a certain degree $d$. The class of ideals having this behavior has been called of {\em sub-homaloidal type}.
The  basic general result about an ideal  $J$ of sub-homaloidal type obtained in \cite[Proposition 4.6]{fat1} was that its initial degree is $s:=(d+1)/2$, where $d$ stands for the degree of the  associated homaloidal type.
Although this was of an enormous avail for the subsequent results, no aspects of the deeper nature of $J$ in the general case were dealt with.

The present paper picks up the discussion in \cite[Section 4]{fat1} by explaining the algebraic and geometric natures of such an ideal $J$,   under the proviso that the three highest virtual multiplicities equal $(s-1)/2$.
Although \cite[Section 4]{fat1} gave a glimpse of geometric flavor bringing in the so-named {\em principal curves} of classical plane Cremona map theory, the finer discussion was not taken up.
Filling the latter gap is a major goal here.
To much of our surprise, the discussion required a thrust into some unexpected subtleties of the homological facet of the ideal in parallel with geometric ideas from fat point theory.
It touches non-trivial aspects of symbolic powers of ideals of fat points, a subject that has recently caught the attention of several authors (see \cite{Bocci}, \cite{ Dumnicki5}, \cite{Guardo}, just to mention a few).
Thus, for example, we are led to lend special attention to the behavior of the {\em saturation degree} of the second and third symbolic powers of a an ideal of fat points of sub-homaloidal type. Unfortunately, the well-known general optimal bounds of this invariant in the class of unmixed ideals in $R=k[x,y,z]$ will not suit our purpose, so we are lead to develop particular techniques to deal with the present special case.
We believe it may be fair to say that the symbolic power theory of a fat ideal of sub-homaloidal type in degree $s$ has its main features already in low order.
This is the case at least under the proviso that the three highest virtual multiplicities equal $(s-1)/2$. 
It is plausible that in general some asymptotic behavior of symbolic powers may intervene, such as the Waldschmidt constant and questions arising from the optimal generation of the symbolic power algebra.

\smallskip

Next is a  description of the sections. 

The first section is a recap of terminology -- for further details we refer to \cite[Introduction]{fat1}.
In addition, it shortly explains the notion of a source inversion factor and its use within the discussion of symbolic powers.
The notion of sub-homaloidal type and its relation to the geometric nature of the corresponding second symbolic power type is revisited.
Although the latter matter has been partly discussed in \cite{fat1}, this time around we list the main open questions regarding the class of fat ideals of sub-homaloidal type.

\smallskip

The second section collects auxiliary results of varied nature, mostly related to the use of arithmetic quadratic transformations.
We start (first subsection) with a brief account of this classical method, adding an algebraic lemmata thereof that clarifies some of the properties in the passage from the given vector of virtual multiplicities to its transform.
This has been useful in deriving other results and by and large we believe it is a very useful tool.
Arithmetic quadratic transformations are of course a venerable topic with a modern view by Nagata. 
Some of it is freely used in the current literature without source quoting or even quoted under different name (e.g., in  \cite{CiMi} it is called {\em Cremona equivalence}).
Of course its geometric usefulness is clear as it tantamount to transforming plane linear systems or divisors by applying a basic plane Cremona transformation.

In the second subsection we apply these transformations to sub-homaloidal types.
In the case of a sub-homaloidal type of the form $(s;(s-1)/2^3,\mu_4,\ldots,\mu_r)$, such that the corresponding homaloidal type is proper, we elaborate on the properties of the transformed type. In particular, we give a complete description of the numerical and the homological nature of the transformed fat ideal $\widetilde{J}$ of the associated fat ideal $J$. Thus, we show, among other things, that the initial degree of $\widetilde{J}$ is $(s+3)/2$ and the linear system in this degree has the expected dimension (which is one plus the initial degree); besides, we prove that   the minimal free resolution of $\widetilde{J}$ is linear.

The third subsection is the bulkiest.
The section brings the diversity of results of the previous sections closest to the main goal of the paper.
This is done on three main fronts: one for the fat ideal $J$ of a (proper) sub-homaloidal type of the form $(s;(s-1)/2^3,\mu_4,\ldots,\mu_r)$, one for the base ideal of the Cremona map derived from $J^{(2)}$, and one for the underlying geometry of the transformed fat ideal $\widetilde{J}$. 
Here we first show that the rational map $\mathfrak{F}$ defined by a minimal set of generators of $\widetilde{J}$ is birational onto the image. In addition, we prove that the initial degree of the second symbolic power $\widetilde{J}^{(2)}$ is $s+2$ and that a vector basis of its linear system at this degree contains $s$ independent source inversion factors of  $\mathfrak{F}$.
We conjecture that actually  $\widetilde{J}^{(2)}$ is generated in degree $s+2$ and that the $s$ source inversion factors form a basis.  

Second we show that $J$ is in fact generated in degree $s$ and describe its minimal resolution. Further a basis of $J_s$ defines a birational map onto the image.
We prove that the image $V$ has degree $s$ in this embedding and, in addition, drawing upon a result of \cite{GGP} as applied to the scheme version of $\widetilde{J}$ we deduce that the $V_s$ is an arithmetically Cohen--Macaulay variety.
As a consequence we derive that the reduction number of the ideal $J$ is at most $2$.
Finally, using all these resources we apply a criterion of Cortadellas--Zarzuela to show that the Rees algebra of the ideal $J$ is Cohen--Macaulay. This implies via a theorem of Aberbach--Huneke--Trung that the minimal defining equations of $V_s$ have degree at most $3$.

Third, we relate the previous results back to the structure of the base ideal of the original Cremona map that motivated them.
Here we obtain the details of the minimal free resolution of this ideal.
Although the shape of the resolution had been dealt with in \cite[Propostion 4.12]{fat1}, the method employed here tells how birational invariants come into the picture.
Namely, consider the inverse map to the birational map defined by $J_s$.
This map has three distinctive representatives yielding three respective source inversion factors $D_1,D_2,D_3$ of degree $2s-1$ each. Then the base ideal is generated by $D_1,D_2,D_3$. The rest of the resolution is derived as in \cite{fat1}.

\smallskip

The last subsection of the section hinges on the so-called Bordiga--White parametrizations studied in \cite{GGP}. The main result consists in showing how a fat ideal of sub-homaloidal type can often produce such parameterizations by applying a quadratic transformation.

\smallskip

As a pointer, the following are the main results of the paper: Proposition~\ref{tilde-generation}, Theorem~\ref{birr_tildeJ},  Theorem~\ref{first_main_theorem}, Theorem~\ref{second_main_theorem} and Proposition~\ref{Bordiga_main}.

\smallskip  

The Appendix at the end gives a few worked-out examples.
 It may essentially be read independently of the rest of the paper and bears no particularly technical tools from ideal theory.
 Some of these are intended to show how one can derive proper sub-homaloidal types  of sufficiently low virtual multiplicities by hand calculation.
 For a numerous list of sub-homaloidal types via computation one can use the routine written by A. Doria in {\em Singular}. A more elaborated routine  returning  sub-homaloidal types whose doubled is proper is also feasible by a suitable implementation of the classical Hudson's test for properness (see \cite[Section 5]{alberich}).
 
 \smallskip
 
 We assume throughout that the ground field has characteristic zero.

\section{Preliminaries}\label{recap}

\subsection{Basic terminology}\label{basic-term}

In this section we briefly recapitulate the basic material to be used in the sequel.
For additional details, see \cite{alberich},  \cite[Introduction]{HS}, \cite[Section 1]{fat1}.

\smallskip

By a {\em linear system} of plane curves of degree $d$ we mean a $k$-vector subspace $L_d$ of the vector space of forms of degree $d$ in the standard graded polynomial ring $R: =k[x,y,z]$.
A linear system $L_d$ defines a rational map $\mathfrak{F}=\mathfrak{F}_d$ with source $\pp^2$ and target $\pp^r$, where $r+1$ is the vector space dimension of $L_d$.

When $\mathfrak{F}$ is birational onto its image, one has important forms attached to the inverse map.
Quite generally, let $\mathfrak{F}:\pp^n\dasharrow \pp^m$ denote a birational map of $\pp^n$ onto the image defined by forms $(f_0:\cdots :f_m)$ of the same degree and
let $(f'_0:\cdots: f'_n)$ denote forms defining the inverse map.
Set $R=k[x_0,\ldots, x_n]$ for the homogeneous coordinate ring of $\pp^n$.
Then there is a uniquely defined
form $D\in R$ such that $f'_i(f_0,\ldots,f_m)=x_iD$, for every $i=0,\ldots,n$.
We call $D$ the {\em {\rm (}source{\rm )} inversion factor} of $\mathfrak{F}$ associated to the chosen representative $(f'_0:\cdots: f'_n)$ of the inverse map.
By  (\cite[Section 2]{AHA}), any such representative, as a vector, can be taken to be
the transpose of a minimal generator of the syzygy module of the so-named weak Jacobian dual matrix.

Note that if the common degree of the $f$'s is $d$ and that of the $f'$'s is $d'$ then the degree of $D$ is $dd'-1.$

The following property of the inversion factor will play a  role in the sequel.
For convenience, we call the ideal of $R$ generated by the defining forms of $\mathfrak{F}$ the {\em base ideal} of this map.

\begin{Proposition}{\rm \cite[Proposition 1.4]{Zaron}}\label{inversionfactor_is_symbolic}
	Let $\mathfrak{F}:\pp^n\dasharrow \pp^m$ be a rational map of degree $d\geq 2$ and let $I\subset R=k[\xx]$ denote its base ideal. Assume that $\mathfrak{F}$ is birational onto the image and let $D\subset R$ denote the source inversion factor relative to a given representative of the inverse map.
	
	Suppose that ${\rm depth}(R/I)>0$. Then:
	\begin{enumerate}
		\item[{\rm (a)}] $D$ is an element of the symbolic power $I^{(d')}$,
		where $d'$ is the degree of the coordinates of the representative.
		In particular, $I^{(d')}\neq I^{d'}$.
		\item[{\rm (b)}] If, moreover, $I^{(\ell)}=I^{\ell}$, $\ell\leq d'-1$, then $D$ is a essential symbolic element of order $d'$.
		In addition,  if $I^{(d')}$ is generated in standard degree $\geq dd'-1$,  then $D$ is a minimal homogeneous generator thereof.
	\end{enumerate}
\end{Proposition}

For a homogeneous ideal $J\subset R$, $h_{R/J}(t):=\dim_k(R_t/J_t)\, (t\geq 0)$ denotes the  Hilbert function of $R/J$.
The following information was gathered in \cite[Section 1.2]{fat1} and is well-known -- here we chose to restate it as a lemma for referencing convenience.

\begin{Lemma}\label{general_unmixed}
Let $J\subsetneq R:=k[x,y,z]$ denote an unmixed homogeneous ideal of codimension $2$.
Then:
\begin{itemize}
\item[{\rm (i)}] $R/J$ is Cohen--Macaulay with minimal graded $R$-resolution of $J$ has the form
\begin{equation}\label{resolution}
0\rar F_1\stackrel{\phi}{\lar} F_0\lar J\rar 0,
\end{equation}
where $\phi$ is a homogeneous map and $F_0,F_1$ are free graded modules.
\item[{\rm (ii)}]  The Hilbert polynomial of $R/J$ is of degree $0$ and coincides with the multiplicity $e(R/J)$.
\item[{\rm (iii)}]  $h_{R/J}(t)$ is strictly increasing till it reaches its maximum $e(R/J)$ and stabilizes thereon.
\item[{\rm (iv)}]  The least $t$ such that $h_{R/J}(t)=e(R/J)$ is the {\em regularity index} of $R/J$
\item[{\rm (v)}]  Set, moreover, $F_0=\oplus_{t>0} R(-t)^{n_t}$, so that the minimal number of generators of $J$ is $\sum_{t>0} n_t$.
Then $n_t=\dim_k{\rm coker}(R_1\otimes_k J_{t-1}\lar J_{t})=\dim_k(J_{t})-\dim_k (R_1J_{t-1})$, where $R_1J_{t-1}$ is short for the vector $k$-subspace of $J_{t}$ spanned by the products of $x,y,z$ by the elements of a basis of $J_{t-1}$ and $R_1\otimes_k J_{t-1}\lar J_{t}$ denote the natural multiplication map.
\end{itemize}
\end{Lemma}

Now we go into more specific behavior regarding an ideal $J:=I(\mathbf{p};\boldsymbol\mu)$ of fat points.
Any such ideal is a perfect ideal, i.e., the ring $R/I(\mathbf{p};\boldsymbol\mu)$ is Cohen--Macaulay.
In addition, the degree of the corresponding scheme  depends only on the number of
points and their appended multiplicities:
\begin{equation}\label{multplicity_of_fat}
e(R/I(\mathbf{p};\boldsymbol\mu))=\sum_{j=1}^n\, \frac{\mu_j(\mu_j+1)}{2}.
\end{equation}

For convenience, we switch to the Hilbert function of the ideal $J$.
Here, for any $t\geq 0$, one has
$$h_J(t):=\dim_k J_t\geq\max\left\{0, {{t+2}\choose {2}} -\sum_{j=1}^n\, \frac{\mu_j(\mu_j+1)}{2}\right\}.$$
When the Hilbert function $h_J(t)$
achieves this lower bound in every degree $t$ one says that it is {\em minimal} -- this is the same as the more common designation in the literature saying that $R/J$  has {\em maximal} Hilbert function. 
A relevant fact is that, due to the behavior in Lemma~\ref{general_unmixed} (iii), $J$ will have minimal Hilbert function if one has
$$h_J(t)={{t+2}\choose {2}} -\sum_{j=1}^n\, \frac{\mu_j(\mu_j+1)}{2},$$
for $t={\rm indeg}(J)$.

As it turns (see Lemma~\ref{general_unmixed} (iv)), the (Castelnuovo--Mumford) regularity of $J$ is the least $t\geq 0$ such that 
$$\dim_k h_J(t-1)= {{t+1}\choose {2}} -\sum_{j=1}^n\, \frac{\mu_j(\mu_j+1)}{2}.$$
It follows that if $J$ has minimal Hilbert function then the regularity is at most ${\rm indeg}(J)+1$, where ${\rm indeg}(J)$ denotes the initial degree of $J$.
This implies by via of \cite{DGM} that, in the notation of Lemma~\ref{general_unmixed} (v), $n_t=0$ except for $t={\rm indeg}(J)$ and possibly $t={\rm indeg}(J)+1$. 

A major question asks when the multiplication map $\mathfrak{m}_t:R_1\otimes_k J_{t-1}\lar J_{t}$ in Lemma~\ref{general_unmixed} (v) has maximal rank. This question has much interest regardless of the assumption that $J$ has minimal Hilbert function, and yet, in the presence of the latter, the crucial multiplication map is the one with $t={\rm indeg}(J)$.
In this case saying that $\mathfrak{m}_{{\rm indeg}(J)}$ is injective or surjective means, respectively, that  $J$ admits no nonzero linear syzygies or that $J$ is generated in one single degree.

The philosophy developed by B. Harbourne in several papers (\cite{HHF}, \cite{H1}, \cite{H2}) is that, in order to obtain information about the  maximal rank property of $\mathfrak{m}_t$ when $J$ has minimal Hilbert function, one needs some elbow room.
The recurrent idea behind this philosophy is to introduce the following additional ideals of fat points:
\begin{equation}\label{j-minus}
\boldsymbol{\mu}^-:=\{\mu_1-1,\mu_2,\ldots, \mu_r\},\: J^-:=I(\mathbf{p}, \boldsymbol{\mu}^-),
\end{equation}
if $\mu_1\geq 1$,
and
\begin{equation}\label{j-plus}
\boldsymbol{\mu}^+:=\{\mu_1+1,\mu_2,\ldots, \mu_r\},\: J^+:=I(\mathbf{p}, \boldsymbol{\mu}^+).
\end{equation}

Among others, Harbourne proves the following important result:

\begin{Lemma}\label{Harbourne_gadgets}{\rm \cite[Lemma 2.4 (d)]{HHF}}
With the above notation, assume that:
\begin{enumerate}
	\item[{\rm (a)}] All three ideals $J$, $J^-$ and $J^+$ have minimal Hilbert function. 
	\item[{\rm (b)}] $h_{J^-}({\rm indeg}(J)-1)>0$ and $h_{J^+}({\rm indeg}(J))>0$.
\end{enumerate}
Then the multiplication map $\mathfrak{m}_{{\rm indeg}(J)}$ is surjective{\rm ;} in particular, $J$ is generated in one single degree.
\end{Lemma}

\subsection{Symbolic square and sub-homaloidal fat ideals}

The following notion was introduced in \cite[Section 4.2]{fat1}.

\begin{Definition}\rm
Let $\boldsymbol\mu=\{\mu_1,\ldots,\mu_r\}$ be a set of nonnegative integers  satisfying the following condition: there exists an integer $s\geq 2$ such that
\begin{equation}\label{eqs-subhomaloidal}
\sum_{i=1}^r \mu_i=3(s-1) \quad {\rm and}\quad \sum_{i=1}^r \mu_i^2=s(s-1).
\end{equation}
Note that the first of the above relations implies that $s$ is uniquely determined.
We will say that $\boldsymbol\mu$ is a {\em  sub-homaloidal multiplicity set in degree $s$}.
\end{Definition}

The basic result relating this notion to the one of a homaloidal multiplicity set is the following:

\begin{Lemma}\label{doubling}
Let $\boldsymbol\mu=\{\mu_1,\ldots,\mu_r\}$ stand for a set of nonnegative integers and let $s\geq 2$ denote an integer.
The following conditions are equivalent:
\begin{enumerate}
\item[{\rm (i)}] $\boldsymbol\mu$ is a sub-homaloidal multiplicity set in degree $s$
\item[{\rm (ii)}] $(2s-1;2\boldsymbol\mu)$ is a homaloidal type.
\end{enumerate}
In addition, fixing a set of points $\mathbf{p}$ with same cardinality as $\boldsymbol\mu$, then $I(\mathbf{p}, 2\boldsymbol\mu)=I(\mathbf{p}, \boldsymbol\mu)^{(2)}$ {\rm (}second symbolic power{\rm )}.
\end{Lemma}

When $\boldsymbol\mu$ is a sub-homaloidal multiplicity set in degree $s$ it  will often be convenient to call the vector $(s;\boldsymbol\mu)$ a {\em sub-homaloidal type}.
Accordingly, we call the ideal $I(\mathbf{p}, \boldsymbol\mu)$ an {\em ideal of sub-homaloidal type}.

In our judgment, the main questions related to this notion are the following;
\begin{Question}\label{main_queries}\rm
	Let $(s, \boldsymbol\mu)$ as above denote a sub-homaloidal type in degree $s$ such that $(2s-1;2\boldsymbol\mu)$ is a {\em proper} type  (i.e., a homaloidal type for which there exist Cremona maps with this type). Fix a set of general points $\mathbf{p}$ with same cardinality as $\boldsymbol\mu$ and let $J:=I(\mathbf{p}, \boldsymbol\mu)$.
	\begin{enumerate}
		\item[{\rm (a)}]  Does $J$ have minimal Hilbert function?
		\item[{\rm (b)}]  Is $J$ generated in degree $s$?
		\item[{\rm (c)}]  Does the linear system $J_s$ define a birational map of $\pp^2$ onto the image?
		\item[{\rm (d)}]  Assuming that (c) holds, is the image of the map ideal theoretically defined by quadrics and cubics?
	\end{enumerate}
\end{Question}
	
Note that having (a) plus knowing that the multiplication map of order $s$ is surjective gives (b).
Question (c) is of course more delicate but, provided (b) takes place,  we expect it to hold in a variety of many situations which allow a certain control of the linear syzygies of $J$.
Finally (d) is a challenging puzzle related, so we see at this stage, to whether the base ideal of a Cremona map with type $(2s-1;2\boldsymbol\mu)$ and general proper base points is saturated or not.

Let us assume once for all that $\mu_1\geq\mu_2\geq\cdots\geq\mu_r$.
It is our view that for a sub-homaloidal multiplicity set in degree $s$ the condition $\mu_i\leq (s-1)/2$, for every $i\geq 1$, is a fundamental virtual barrier.
In the present setup, if $\mu_1>\mu_2$ the inequality $\mu_i\leq (s-1)/2$ for $i\geq 2$ will hold no matter what provided $(2s-1;2\boldsymbol\mu)$ is besides proper.

At the other end, the following question has been stated in \cite[Question 4.1]{fat1}:
\begin{Question}\label{Q1}\rm
Let $(d;\mu_1\geq\cdots\geq\mu_r)$ denote a proper homaloidal type and let $F$ stand for a Cremona map with this type and general base points. If the base ideal of $F$ is not saturated then $\mu_1\leq \left \lfloor{d/2}\right \rfloor$.
\end{Question}


\section{Arithmetic and geometry of subhomaloidal types}

\subsection{Arithmetic quadratic transformations}\label{AQT}

In this short piece we record a few features of this theory needed in the sequel.
Arithmetic quadratic transformations were essentially introduced in the work of Hudson, Coble and Semple--Roth. However, the first systematic treatment is due to Nagata in the monumental \cite{Nagata1} and \cite{Nagata2}.
Here we follow closely the excellent expos\'e of \cite[Section 5.2]{alberich}.

Let $S:=(d;\boldsymbol\nu)$ be a vector in $\mathbb{N}\times \mathbb{N}^r=\mathbb{N}^{r+1}$, where $d>0$ and $\boldsymbol\nu=(\nu_1, \ldots\nu_r)$ -- to fix ideas think of $d$ as a ``degree'' and $\boldsymbol\nu$ as ``virtual multiplicities''.
We will often refer to $\boldsymbol\nu$ as a set rather than a vector, a minor abuse without serious consequences.
In addition we may often call a vector such as $S$ a {\em type} as a reminiscence of the classical terminology in the case $S$ is a homaloidal type.

Fix three indices $j,k,l$ among $\{1,\ldots,r\}$.
The {\em arithmetic quadratic transformation} based on the triple $j,k,l$ is the map
$\mathcal{Q}_{j,k,l}: \mathbb{N}^{r+1}\lar \mathbb{N}^{r+1}$ such that $\mathcal{Q}_{j,k,l}(S)=\widetilde{S}$ with $\widetilde{S}=(\widetilde{d};\widetilde{\boldsymbol\nu})$, where
	$$\left\{
	\begin{array}{l}
	\widetilde{d}:=2d-\nu_j+\nu_k+\nu_l\\
	\widetilde{\nu_j}=d-\nu_k-\nu_l\\
	\widetilde{\nu_k}=d-\nu_j-\nu_l\\
	\widetilde{\nu_l}=d-\nu_j-\nu_k\\
	\widetilde{\nu_i}=\nu_i, \,i\notin\{j,k,l\}.
	\end{array}
	\right.
	$$
Note an important feature of this map: if $\nu_j\geq \nu_k\geq\nu_l$ then $\widetilde{\nu_j}\geq \widetilde{\nu_k}\geq \widetilde{\nu_l}$.

Whenever the indices $j,k,l$ are fixed in the discussion, we denote the transformation simply by $\mathcal{Q}$.

Let ${\bf p}=\{p_1,\ldots,p_r\}$ denote a set of $r$ points in $\pp^2$.
Fixing indices $j,k,l$ as above, we write $\widetilde{\bf p}:=\{\widetilde{p_1},\ldots, \widetilde{p_r}\}$ where 
$$\widetilde{p_i}=\left\{
\begin{array}{ll}
p_i & \mbox{if $i\in\{j,k,l\}$}\\
\mathfrak{Q}(p_i) & \mbox{if $i\notin\{j,k,l\}$}
\end{array}
\right.
$$
with $\mathfrak{Q}$ denoting the plane quadratic Cremona map based on the three points $p_j,p_k,p_l$.

Finally, we let as before $J:=I(\mathbf{p};\boldsymbol\nu)$ and $\widetilde{J}:=I(\widetilde{\mathbf{p}};\widetilde{\boldsymbol\nu})$ denote the respective ideals of fat points in $R=k[x,y,z]$.

Note that for each $d>0$ the integer vector $(d;\boldsymbol\nu)$ is closely associated with the vector space $J_d$ -- i.e., the linear system of forms of degree $d$ passing through each point $p_i$ with multiplicity at least $\nu_i$.
To go beyond this simple association we note that  $J_d$ defines a rational map $\mathcal{F}:\pp^2\dasharrow \pp^N$, for some $N$, and the composite $\mathcal{F}\circ \mathfrak{Q}$ can easily be obtained in algebraic terms.

The next lemma explains how it connects with the rational map defined by $\widetilde{J}_{\widetilde{d}}$.
Parts of the result seem to be taken for granted in the literature, but since we could not find a precise reference, we give a proof.

Note that if  $p_1=(1:0:0), p_2=(0:1:0), p_3=(0:0:1)$, then $\mathfrak{Q}=(yz:xz:yz)$.

\begin{Lemma}\label{transfquadpreserva} 
Assume that $\nu_1\geq \ldots\geq\nu_r$ and let $j=1,k=2.l=3$.
Let ${\bf p}=\{p_1,\ldots,p_r\}$ be points in $\pp^2$, the first three of which are the coordinate points and let $J, \widetilde{J}$ be defined as above.
Set $\Delta:=x^{\nu_1}y^{\nu_2}z^{\nu_3}$. Then, for any given $d>0$ one has:
\begin{enumerate}
\item[{\rm (a)}] If $f(x,y,z)\in J_d$ then $\Delta$ divides $f(yz:xz:yz).$
\item[{\rm (b)}] If $f(x,y,z)\in J_d$ then $f(yz:xz:yz)/\Delta\in \widetilde{J}_{\widetilde{d}}.$
\item[{\rm (c)}] For any integer $n\geq 1$, the map $\eta^{(n)}:f(x,y,z)\mapsto f(yz:xz:yz)/\Delta^n$ induces an isomorphism of the $k$-vector spaces  $(J^n)_{nd}$ and  $(\widetilde{J}^n)_{n\widetilde{d}}.$
\end{enumerate}
\end{Lemma}
\demo
(a) Since $f\in (y,z)^{\nu_1},$ then  $f(yz,xz,xy)\in(x^{{\nu_1}}).$ Similarly, $f(yz,xz,xy)\in(y^{{\nu_2}})$ and $f(yz,xz,xy)\in(z^{{\nu_3}}).$ Thus, $\Delta$ is a factor of $f(yz,xz,xy)$ as claimed.

\smallskip

(b) We have to show that the multiplicity of $f(yz,xz,xy)/\Delta$ at each $p_i$ is at least $\widetilde{\nu}_i.$ Since $f$ is homogeneous of degree $d$,  $f(yx,xz,yz)\in (y,z)^{d}\cap(x,z)^{d}\cap(x,y)^d.$ Thus,

$$e_{p_1}\left(\frac{f(yz,xz,xy)}{\Delta}\right)=e_{p_1}(f(yz,xz,xy))-e_{p_1}(\Delta)\geq d-(\nu_2+\nu_3)=\widetilde{\nu}_1,$$
$$e_{p_2}\left(\frac{f(yz,xz,xy)}{\Delta}\right)=e_{p_2}(f(yz,xz,xy))-e_{p_2}(\Delta)\geq d-(\nu_1+\nu_3)=\widetilde{\nu}_2$$
e
$$e_{p_3}\left(\frac{f(yz,xz,xy)}{\Delta}\right)=e_{p_3}(f(yz,xz,xy))-e_{p_3}(\Delta)\geq d-(\nu_1+\nu_2)=\widetilde{\nu}_3$$

Next let $p=(\alpha_1:\alpha_2:\alpha_3)$ denote any of the remaining points $p_i,$ with $4\leq i\leq r,$ and let $\nu$ denote the respective virtual multiplicity. Since the defining ideal of $p$ is $(\alpha_3x-\alpha_1z,\alpha_3y-\alpha_2z)$ then that of $\mathfrak{Q}(p)$ is $(\alpha_1x-\alpha_3z,\alpha_2y-\alpha_3z).$
Then one must have

$$f=\sum_{i=1}^{\nu}a_{i}(x,y,z)(\alpha_3x-\alpha_1z)^{i}(\alpha_3y-\alpha_2z)^{\nu-i}$$
for certain $a_i\in k[x,y,z].$ 
It follows that

$$f(yz,xz,xy)=\sum_{i=1}^{\nu}a(yz,xz,xy)y^{i}x^{\nu-i}(\alpha_3z-\alpha_1x)^{i}(\alpha_3z-\alpha_2y)^{\nu-i}.$$
Since $(\alpha_3z-\alpha_1x)^{i}(\alpha_3z-\alpha_2y)^{\nu-i}\i$ is the typical term in the binomial expansion of  the power $(\alpha_1x-\alpha_3z,\alpha_2y-\alpha_3z)^{\nu}$, then $f(yz,xz,xy)\in(\alpha_1x-\alpha_3z,\alpha_2y-\alpha_3z)^{\nu},$ which implies that  $e_p(f(yz,xz,xy))\geq\nu.$ 
Therefore,
$$e_{p}\left(\frac{f(yz,xz,xy)}{\Delta}\right)=e_{p}\left(f(yz,xz,xy)\right)-e_{p}\left(\Delta\right)\geq \nu$$
as claimed. 

\smallskip

(c) Set $\widetilde{\Delta}:=x^{\widetilde{\nu}_1}y^{\widetilde{\nu}_2}z^{\widetilde{\nu}_3}.$ and consider the analogous map 
$$\zeta^{(n)}: g(x,y,z)\mapsto g(yz,xz,xy)/\widetilde{\Delta}$$ 
coming the other way. By the same token, it maps
$\widetilde{J}_{\widetilde{d}}$ to $J_d$.
We claim that the two maps are inverse of each other. 

First consider the case $n=1$, setting $\eta^{(1)}=\eta$ and $\zeta^{(1)}=\zeta$ for lighter reading.
One has
\begin{eqnarray}
\zeta\circ\eta(f(x,y,z))&=&\zeta(f(yz,xz,xy)/\Delta)\nonumber\\
&=&\frac{f((xyz)x,(xyz)y,(xyz)z)}{\Delta(yz,xz,xy)\widetilde{\Delta}}\nonumber\\
&=&\frac{x^dy^dz^df(x,y,z)}{(x^{d-\widetilde{\nu}_1}y^{d-\widetilde{\nu}_2}z^{d-\widetilde{\nu}_3})(x^{\widetilde{\nu}_1}y^{\widetilde{\nu}_2}z^{\widetilde{\nu}3})}\nonumber\\
&=&f(x,y,z)\nonumber
\end{eqnarray}
The desired equality in the other direction is the same calculation.

Now take the case of arbitrary $n\geq 1$.
Let $\{f_1,\ldots,f_m\}$ denote a vector basis of $J_d.$ Since any element of $(J^n)_{nd}$ is a $k$-linear combination of the forms $\{f_1^{\alpha_1}\cdots f_m^{\alpha_m}\;|\;\alpha_1+\cdots+\alpha_m=n\},$ it suffices to deal with one such form.
Clearly, 
 $$\eta^{(n)}(f_1^{\alpha_1}\cdots f_m^{\alpha_m})=\eta(f_1)^{\alpha_1}\cdots \eta(f_m)^{\alpha_m}$$
A similar rule works for $\zeta^{(n)}$ as well.
Therefore
$$\zeta^{(n)}(\eta^{(n)}(f_1^{\alpha_1}\cdots f_m^{\alpha_m}))=\zeta(\eta(f_1))^{\alpha_1}\cdots \zeta(\eta(f_m))^{\alpha_m}=(f_1)^{\alpha_1}\cdots (f_m)^{\alpha_m}$$
from the case $n=1$.

The required equality in the other direction works the same way.
\qed

\subsection{$3$-uniform subhomaloidal types} 

Classically, arithmetic quadratic transformations act on integer vectors of $\mathbb{N}^{r+1}$ which are homaloidal or exceptional types.
In this work they will also act on vectors which are subhomaloidal types, as will be seen in this subsection.

We fix the following general notation:

\begin{itemize}
	\item $\boldsymbol\mu=\{\mu_1\geq\cdots \geq\mu_r\}$ denotes a sub-homaloidal multiplicity set in degree $s\geq 3$.
	\item $S:=(s;\boldsymbol\mu)$ denotes the corresponding sub-homaloidal type.
	\item $S^{(2)}:=(2s-1;2\boldsymbol\mu)$ denotes the ``doubled'' homaloidal type.
	\item $\mathcal{Q}$ and $\mathfrak{Q}$ respectively denote, as in the previous subsection, the arithmetic quadratic transformation  based on the thee highest virtual multiplicities $\mu_1\geq\mu_2\geq\mu_3$ and the corresponding uniquely defined plane quadratic Cremona map with proper base points $p_1,p_2,p_3$.
	\item $\widetilde{S}:=(\widetilde{s};\widetilde{\boldsymbol\mu})$ will be the transform of $S$ by  $\mathcal{Q}$, as explained in the previous subsection.
	\item $\mathbf{p}=\{p_1,\ldots,p_r\}\subset \pp^2$ denotes a set of $r$ general points and $$\widetilde{\mathbf{p}}, J=I(\mathbf{p};\boldsymbol\mu), \widetilde{J}:=I(\widetilde{\mathbf{p}};\widetilde{\boldsymbol\mu})$$ will have the same meaning as already explained in the previous section.
\end{itemize}

In such a generality not much can be said, except for the following result:

\begin{Lemma}{\rm (\cite[Proposition 4.6]{fat1}}\label{indeg_of_J}
With the above notation, if the homaloidal type $S^{(2)}$ is proper then $s$ is the initial degree of $J$.
\end{Lemma}

Lemma~\ref{transfquadpreserva} gives a basic interplay between the linear systems $J_d$ and $\widetilde{J}_{\widetilde{s}}$.
But in fact one can go deeper under tighter hypotheses as we next examine.
Namely, we will  assume once for all throughout this subsection the following hypothesis: {\em the homaloidal type $S^{(2)}$ is proper and} $\mu_1=\mu_2=\mu_3= (s-1)/2$ -- a situation we informally designate as $3$-uniform by mimicking the existing terminology.

\begin{Proposition}\label{tilde-generation}
With the above notation, let $S^{(2)}$ denote a proper $3$-uniform homaloidal type. Then:
	\begin{enumerate}
		\item[{\rm (a)}] ${\rm indeg}(\widetilde{J})=(s+3)/2$.
		\item[{\rm (b)}] $\widetilde{J}$ has linear minimal resolution
		\begin{equation}\label{lin-res}
		0\rar R^{(s+3)/2}(-((s+5)/2))\lar R^{(s+5)/2}(-(s+3)/2)\rar \widetilde{J}\rar 0.
		\end{equation}
		In particular, ${\rm reg}(\widetilde{J})=(s+3)/2$.
	\end{enumerate}
\end{Proposition}
\demo (a)
Under the hypothesis on the three highest virtual multiplicities, one has $\tilde{s}=2s-(3/2)(s-1)=(s+3)/2$ and
\begin{equation}\label{type_auxiliary}
\tilde{S}=((s+3)/2; \mu_4,\ldots,\mu_r, 1^3).
\end{equation}

By Lemma~\ref{transfquadpreserva} we know that $\dim J_s=\dim\widetilde{J}_{(s+3)/2}.$

We draw upon parts of the argument in the proof of \cite[Lemma 4.10]{fat1}.
For convenience, recall the main steps.

The focus is on the types
\begin{equation}\label{homaloidal_type}
S^{(2)}=(2s-1;(s-1)^3,2\mu_4,\ldots,2\mu_r)
\end{equation}  
and
\begin{equation}\label{first_exceptional_type}
E:=(s-1; ((s-1)/2)^2,(s-3)/2,\mu_4,\ldots,\mu_r).
\end{equation}  
Applying $\mathcal{Q}$ as usual yields 
\begin{equation}\label{homaloidal_transformed}
\widetilde{S^{(2)}}=(s+1;2\mu_4,\ldots,2\mu_r,1^3)
\end{equation} 
and
\begin{equation}\label{second_exceptional_type} \widetilde{E}=((s+1)/2;\mu_4,\ldots,\mu_r,1^2, 0),
\end{equation} 
respectively.
Since $\widetilde{S^{(2)}}$ is still a proper homaloidal type then $\widetilde{E}$ is a proper exceptional type, as shown in the proof of \cite[Lemma 4.10]{fat1}.
Set 
$$\widehat{J}:=I(p_2)\cap I(p_3)\cap I(\mathfrak{Q}(p_4))^{\mu_4}\cap\ldots\cap I(\mathfrak{Q}(p_r))^{\mu_r}.$$

By a similar token and by the nature of the proper exceptional type, one has ${\rm indeg}(\widehat{J})=(s+1)/2$ and $\widehat{J}_{(s+1)/2}$ has dimension $1$, spanned by the equation of a principal curve.

Now, by construction, one has the inclusion of ideals $\widetilde{J}\subset \widehat{J}$, which induces  inclusions at the level of  corresponding graded parts.
It follows that ${\rm indeg}(\widetilde{J})\geq (s+1)/2$.
Since $\widehat{J}_{(s+1)/2}$ is spanned by the equation of a principal curve and the latter passes in each point with effective multiplicity equal to the virtual one, then $\widetilde{J}_{(s+1)/2}\neq \widehat{J}_{(s+1)/2}$, and hence ${\rm indeg}(\widetilde{J})\geq (s+1)/2 +1=(s+3)/2$.
On the other hand, one has $\dim \widetilde{J}_{(s+3)/2}=\dim J_s$ by Lemma~\ref{transfquadpreserva} and
$$\dim J_s\geq {s+2\choose 2}-\sum_{i=1}^{r} {\mu_i+1\choose 2 }=\frac{s+5}{2}=\frac{s+3}{2}+1,$$
where the rightmost term is the expected dimension of $J_s$.
Setting $N:=(s+3)/2$ for lighter reading, one now has that $\widetilde{J}$ contains at least $N+1$ minimal generators of degree $N$. Since $\widetilde{J}$ is minimally generated by the maximal minors of an $M\times (M-1)$ matrix, if $M\geq N+2$ then every maximal minor of this matrix would have degree at least $M-1 \geq N+1$, which is a contradiction.

Therefore, it follows that $\widetilde{J}$ is minimally generated by a vector basis of $\widetilde{J}_{(s+3)/2}$ and, in particular, it is generated in degree $(s+3)/2$.
As a bonus, we have also shown that $J_s$ has the expected dimension $(s+5)/2$.

\medskip

(b) Set $N:=(s+3)/2$ as in the proof of (a).
By (a),   $\widetilde{J}$ is generated by the maximal minors of an $(N+1)\times N$ matrix which necessarily have linear entries. 
This implies that $\widetilde{J}$ has minimal free resolution of the form
\begin{equation*}
\hspace{3cm} 0\rar R^{N}(-(N+1))\lar R^{N+1}(-N)\rar \widetilde{J}\rar 0. \hspace{4cm} \square
\end{equation*}

Besides granting once for all the assumptions and notation of this subsection,  we will maintain the notation $N:=(s+3)/2$.

\begin{Corollary}\label{tilde-symbolic-powers}
	For any integer $n\geq 1$ one has:
	\begin{enumerate}
		\item[{\rm (a)}] ${\rm reg}(\widetilde{J}^n)\leq nN$
		\item[{\rm (b)}] $(\widetilde{J}^{(n)})_{nN}=(\widetilde{J}^n)_{nN}$.
		\item[{\rm (c)}] ${\rm reg}(\widetilde{J}^{(n)})\leq nN.$
		\item[{\rm (d)}] $\dim_{k} \widetilde{J}_{ns}^{n}=\dim_k \widetilde{J}^{(n)}_{n N}= \frac{s}{2}n^2+\frac{3}{2}n+1. $
	\end{enumerate}		
\end{Corollary}
\demo (a) By \cite[Theorem 1.1]{GGP} one has $${\rm reg}(\widetilde{J}^n)\leq n\cdot{\rm reg}(\widetilde{J}).$$ 
Since ${\rm reg}(\widetilde{J})=N$ from Proposition~\ref{tilde-generation} (b), we are done.

\smallskip

(b) 
For a homogeneous ideal $I\subset k[x,y,z],$  let $\mathfrak{n}(I)$ denote its {\em saturation degree}, i.e., the least positive integer $s$  such that $I_t=I^{\rm sat}_{t}$ for every $t\geq s.$ 
One has
$${\rm reg}(I)=\max\{\mathfrak{n}(I),{\rm reg}(I^{\rm sat})\}.$$ 
The result follows immediately from this and item (a).

\smallskip

(c) The proof is the same as in item (b).

\smallskip

(d) Using the hypothesis that $\boldsymbol\mu$ is subhomaloidal and that $\mu_1=\mu_2=\mu_3=(s-1)/2$, a straightforward calculation yields the following relations
	\begin{equation}
	\sum_{i=1}^r (n\widetilde{\mu_i})= 3n(N-1)\quad\mbox{e}\quad \sum_{i=1}^r (n\widetilde{\mu_i})^2=n^2(N-2)N+3n^2,
	\end{equation}
for every integer $n\leq 1$.

From this, $e(R/\widetilde{J}^{(n)})=\frac{n^2(N-2)N+3n^2+3n(N-1)}{2}.$

On the other hand, quite generally for an unmixed ideal $I\subset R$, if $\rho:={\rm reg}(I)\leq D$, for some integer $D$, then $\dim I_{\rho-1}={{\rho +1}\choose {2}}-e(R/I)$, hence the Hilbert function of $R/I$ attains its maximum at $\rho-1$ and then also at $D\geq \rho$. Consequently, the Hilbert function of $I$ at
$D$ has the expected value ${{D+2}\choose {2}}-e(R/I)$.

Applying this in the present case with reg$(\widetilde{J}^{(n)})\leq nN$ yields that $\widetilde{J}^{(n)}_{n N}$ has the expected dimension, i.e,
\begin{eqnarray*}
\dim_k \widetilde{J}^{(n)}_{n N}&=&{nN+2\choose 2}-e(R/\widetilde{J}^{(n)})\\
&=& \frac{(nN+2)(nN+1)}{2}-\frac{n^2(N-2)N+3n^2+3n(N-1)}{2}\\
&=&\frac{2Nn^2-3n^2+3n+2}{2}\\
&=&\frac{3N-3}{2}n^2+\frac{3}{2}n+1\\
&=&\frac{s}{2}n^2+\frac{3}{2}n+1
\end{eqnarray*}
\qed

\medskip

\begin{Proposition}\label{dimensions}
		For any integer $n\geq 1$ one has:
		\begin{enumerate}
			\item[{\rm (a)}] The arithmetic transform of the type $(ns;n\boldsymbol{\mu}$) by $\mathcal{Q}$ is the type $(nN;n \widetilde{\boldsymbol\mu})$
			\item[{\rm (b)}] $J^{(n)}_{ns}=J^n_{ns}$.
			\item[{\rm(c)}] $\dim_{k} J_{ns}^{n}=\frac{s}{2}n^2+\frac{3}{2}n+1. $
		\end{enumerate}		
\end{Proposition}
\demo
(a) This is an obvious calculation from the definitions:
$$2ns-n\mu_1-n\mu_2-n\mu_3=n(2s-\mu_1-\mu_2-\mu_3)=nN,$$
$$ns-n\mu_2-n\mu_3=n(s-\mu_2-\mu_3)=n\widetilde{\mu}_1,$$
$$ns-n\mu_1-n\mu_3=n(s-\mu_1-\mu_3)=n\widetilde{\mu}_2,$$
$$ns-n\mu_1-n\mu_2=n(s-\mu_1-\mu_2)=n\widetilde{\mu}_3$$
and
$$n\mu_i=n\widetilde{\mu}_i,\:\mbox{for}\; i\geq 4.$$

\smallskip

(b) Since $J_{ns}^{n}\subset J_{ns}^{(n)}$, it suffices to show that these are $k$-vector spaces of same dimension.
One has
\begin{eqnarray}
\dim_k (J^{n})_{ns}
&=&\dim_k (\widetilde{J}^n)_{nN}\nonumber\\
&=&\dim_k (\widetilde{J}^{(n)})_{nN}\nonumber\\
&=& \dim_k (J^{(n)})_{ns},\nonumber
\end{eqnarray}
where the first equality comes from Lemma~\ref{transfquadpreserva}~(d) and the last two come from Corollary~\ref{tilde-symbolic-powers}~(b).

\smallskip

(c) This is a straightforward consequence of Corollary~\ref{tilde-symbolic-powers} (d) and from the equality $\dim_k (J^{n})_{ns}
=\dim_k (\widetilde{J}^n)_{nN}.$
\qed

\begin{Proposition}\label{image_has_degree_s}
Let $W\subset \pp^N$ denote the image of the rational map defined by the linear system $J_s$. Then the Hilbert function of $W$ is given by
$H_W(n)=\frac{s}{2}n^2+\frac{3}{2}n+1,\, n\geq0$.
Thus, the Hilbert polynomial of $W$ coincides with its Hilbert function from $n\geq 1;$ in particular, the degree of $W$ is $s$.
\end{Proposition}
\demo
The homogeneous coordinate ring $C$ of $W$ in its embedding is isomorphic by a degree preserving isomorphism to the homogeneous $k$-subalgebra $k[f_0t,\ldots,f_Nt]\subset R[t]$ generated in degree $1$, where $J_s$ is panned by $\{f_0,\ldots,f_N\}$.
Therefore, for any $n\geq 1$ one obtains $k$-vector space isomorphisms
\begin{eqnarray*}
C_n &\simeq & k[f_0t,\ldots,f_Nt]_n\simeq \sum_{|\alpha|=n}\, k (\ff t)^{\alpha}= t^n\,\sum_{|\alpha|=n}\, k (\ff)^{\alpha}\\
&\simeq &\sum_{|\alpha|=n}\, k (\ff)^{\alpha}\simeq  k[f_0,\ldots,f_N]_n\simeq (J^n)_{ns}.
\end{eqnarray*}
Thus, the result follows from Proposition~\ref{dimensions} (c).
\qed

\subsection{The plus/minus Harbourne criterion}

This short part shows that the subhomaloidal type under consideration fulfills Harbourne's criterion in Lemma~\ref{Harbourne_gadgets} as exposed in Subsection~\ref{basic-term}.

\begin{Proposition}\label{Harbourne-holds}
	Notation as in the beginning of the subsection. Let $J^-$ and $J^+$ be as in {\rm (\ref{j-minus})}  and {\rm (\ref{j-plus})}, respectively.
	Then hypotheses {\rm (a)} and {\rm (b)} of {\rm Lemma~\ref{Harbourne_gadgets}} hold for $J, J^-,J^+$.
\end{Proposition}
\demo
Thus, we are to show that the three fat ideals $J, J^-, J^+$ have minimal Hilbert function and that $h_{J^-}({\rm indeg}(J)-1)>0$ and $h_{J^+}({\rm indeg}(J))>0$.
As explained in the first section, in order to show that the Hilbert function of a fat ideal $I:=I(\mathbf{p},\boldsymbol{\mu})$ is minimal it suffices to show that
$$h_I(\alpha)={\alpha+2\choose 2}-\sum_{i=1}^r{\mu_i+1\choose 2},$$
where $\alpha:={\rm indeg}(J)$.

For $J$, it follows from the proof of Proposition~\ref{tilde-generation} (a) that $h_J(s)$ has the expected dimension, hence $J$ has minimal Hilbert function.

As for $J^-$, note that it coincides with the fat ideal on the virtual multiplicities of the type (\ref{first_exceptional_type}) based on the originally given points $p_i, 1\leq i\leq r$.
By the same token as argued about (\ref{second_exceptional_type}), the former is a proper exceptional type, in particular one has
${\rm indeg}(J^-)=s-1$ and $\dim h_{J^-}(s-1)=1$.
On the other hand, a direct calculation gives
$${s+1\choose 2}-{\mu_1\choose 2}-\sum_{i=2}^r{\mu_i+1\choose 2}=1.$$
Since the left-hand side of the above equality is at most $\dim h_{J^-}(s-1)$, then we get equalities throughout.
This shows both that $J^-$ has minimal Hilbert function and that its value on ${\rm indeg}(J)-1$ is positive, as required in item (b) of Lemma~\ref{Harbourne_gadgets}.

We are left with $J^+$.
The easy piece is the equality
$${s+2\choose 2}-{\mu_1+2\choose 2}-\sum_{i=2}^r{\mu_i+1\choose 2}=2,$$
again a direct calculation using the hypothesis that $S^{(2)}=(2s-1;2\boldsymbol{\mu})$  is a proper homaloidal type.

To close the cycle as above we now show that $h_{J^+}(s)=2$ and ${\rm indeg}(J^+)=s$.
We refer once more to the proper exceptional type $\widetilde{E}=((s+1)/2;\mu_4,\ldots,\mu_r,1^2, 0)$ in (\ref{second_exceptional_type}). Applying iterated arithmetic quadratic transformations based on the three highest virtual multiplicities we eventually reach the type $(1;1^2,0^{r-2})$ (see \cite[Proof of Lemma 4.10]{fat1}).

Now consider the type $(s;\boldsymbol{\mu ^+})$, where $\boldsymbol{\mu ^+}$ are the virtual multiplicities of $J^+$.
Applying to this type the arithmetic quadratic transformation based on the three highest virtual multiplicities gives the type 
$$((s+1)/2;\mu_4,\ldots,\mu_r,1, 0^2).$$
Then, as above, applying iterated arithmetic quadratic transformations based on the three highest virtual multiplicities we eventually reach the type $(1;1,0^{r-1})$.
The fat ideal on the virtual multiplicities of the latter type, based on the current set of transformed point, clearly has vector dimension $2$.
Since the iterated procedure preserves the vector dimension of the graded part of the fat ideal in the degree of the type, we get $h_{J^+}(s)=\dim_k(J^+)_s=2$, as desired.

To conclude, we now have ${\rm indeg}(J^+)\leq s$.
If $h_{J^+}(s-1)>0$ let $0\neq f\in J^+_{s-1}$.
Then $\{xf,yf,zf\}\subset J^+_s$ which is absurd for reasons of vector space dimension.
Therefore, ${\rm indeg}(J^+)= s$, as desired.
\qed

\smallskip

As a consequence of Lemma~\ref{Harbourne_gadgets} we get:

\begin{Corollary}\label{equi-generation}
If $J$ is a proper subhomaloidal $3$-uniform fat ideal then the multiplication map $\mathfrak{m}: R_1\otimes_k J_s\rar J_{s+1}$ is surjective$;$ in particular,  $J$ is generated in degree $s$.
\end{Corollary}

 \section{Main results}

\subsection{The structure of the second symbolic power}
	
In the previous part we emphasized the relation between the ordinary and the symbolic powers of the subhomaloidal fat ideal $J$ and its transform $\widetilde{J}$ via an arithmetic quadratic transformation.
In this part we focus mainly on $\widetilde{J}$, giving additional details about its second symbolic power. The obtained structure will have impact on the search of a minimal reduction of $J$ and the corresponding reduction number, a matter of crucial use in the proof of the main theorem in the next section.

\medskip

In the sequel we keep the notation of the previous subsection.

We will make use of the following abstract result.

\begin{Lemma}\label{4guys}
	Let $f,g\in R=k[x,y,z]$ denote forms of the same degree $\geq 2$ with no proper common factor. Then the linear system $L=\{xf,yf,xg,zg\}$ defines a birational map of $\pp^2$ onto the image.
\end{Lemma}
\demo
The argument uses the general birationality criterion established in \cite{AHA}.
For this one considers the Rees algebra of the ideal $(L)$ and let $t,u,v,w$ denote presentation variables of this $R$-algebra.
The linear syzygies $y.xf-x.yf=0$ and $z.xg-x.zg=0$ give forms $ys-xt, zu-xv$ in the presentation ideal of the Rees algebra.
Taking their derivatives with respect to $x,y,z$ yields the matrix
$$\left(\begin{array}{ccc}
-t&s&0\\
-v&0&u
\end{array}\right)
$$
Now, by definition this is a submatrix of the so-called Jacobian dual matrix introduced in \cite{AHA}.
By \cite[Theorem 2.18]{AHA}, the map is birational if and only if the rank of the latter matrix is $2$ modulo the homogeneous defining ideal $P$ of
the image. Since the $2$-minors are monomials  while $P$ is a prime ideal, we will be through if we show that $P$ admits no $1$-forms, i.e., if  $xf,yf,xg,zg$ are $k$-linearly independent.
But the latter is clear by a simple calculation using that $\gcd(f,g)=1$ and that $\deg(f)=\deg(g)\geq 2$.
\qed

\medskip

\begin{Theorem}\label{birr_tildeJ} Let $\mathfrak{F}:\pp^2\dasharrow \pp^{\frac{s+3}{2}}$ denote the rational map defined by a set of minimal generators of $\widetilde{J}.$ Then:
\begin{enumerate}
\item[{\rm (a)}] $\mathfrak{F}$ is birational onto the image.
\item[{\rm (b)}] $\dim _k((\widetilde{J}^{(2)})_{s+2})\geq s$ and the initial degree of $\widetilde{J}^{(2)}$ is $s+2$.
\item[{\rm (c)}] A vector basis of $(\widetilde{J}^{(2)})_{s+2}$ contains $s$ source inversion factors of $\mathfrak{F}$.
\end{enumerate}
\end{Theorem}
\demo (a) 
Since $\widetilde{J}=(\widetilde{J}_{(s+3)/2})$ by Proposition~\ref{tilde-generation}, then the map is defined by a $k$-basis of the linear system $\widetilde{J}_{(s+3)/2}$.
Still by the proof of Proposition~\ref{tilde-generation}, the vector $\widetilde{E}=((s+1)/2;\mu_4,\ldots,\mu_r,1^2, 0)$ in (\ref{second_exceptional_type}) is an exceptional type for either version of the virtual multiplicities
$$\boldsymbol\mu_1=(0,1,1,\mu_1,\ldots,\mu_r)  \; {\rm and} \;\boldsymbol\mu_2=(1,0,1,\mu_1,\ldots,\mu_r).$$
Thus, the linear systems $I(\widetilde{\bf p};\boldsymbol\mu_1)_{(s+1)/2}$ and $I(\widetilde{\bf p};\boldsymbol\mu_2)_{(s+1)/2}$  are spanned by the equations $f$ and $g$ of principal curves, respectively.
Since $f,g$ are irreducible and necessarily distinct then by Lemma~\ref{4guys}, $\{xf,$ $yf,$ $xg,zg\}$ defines a birational map of $\pp^2$  onto the image.
Clearly,  $\{xf, yf, xg,zg\}\subset\widetilde{J}_{(s+3)/2}.$ 
Therefore, $\widetilde{J}_{(s+3)/2}$ also defines a birational map of $\pp^2$ onto the image.

\medskip

(b) Recall that $\widetilde{J}^{(2)}$ is the fat ideal with same base points as and doubled virtual multiplicities.
Therefore, it suffices to show that the expected dimension of the linear system $(\widetilde{J}^{(2)})_{s+2}$ is $s$.
Having in mind that 
$$\widetilde{\nu}_i=\left\{
\begin{array}{ll}
1 & \mbox{if $i\in\{1,2,3\}$}\\
\nu_i & \mbox{if $i\geq 4$}
\end{array}
\right.,
$$
and drawing upon the subhomaloidal equations (\ref{eqs-subhomaloidal}) of $J$, one obtains

\begin{eqnarray*} h_{s+2}(\widetilde{J}^{(2)})&=&{{s+2}\choose {2}}-\sum_i \frac{2\widetilde{\nu}_i(2\widetilde{\nu}_i+1)}{2}=9+\sum_{i\geq 4}\nu_i+2\sum_{i\geq 4}\nu_i\\
	&=&{{s+2}\choose {2}}- 9 -(3s(s-1)-\frac{3}{2}(s-1)) -2(s(s-1)-\frac{3}{4}(s-1)^2)\\
	&=&{{s+2}\choose {2}}-\frac{s^2+5s+12}{2}=s.	
\end{eqnarray*}

\smallskip

In particular, $(\widetilde{J}^{(2)})_{s+2}\neq 0$. Thus, for the assertion on the initial degree it suffices to show that $(\widetilde{J}^{(2)})_{s+1}= 0$.
For this we go back to the level of $J$ and consider the transformed homaloidal type $\widetilde{S^{(2)}}=(s+1;2\mu_4,\ldots,2\mu_r,1^3)$ in (\ref{homaloidal_transformed}).
By the proof of \cite[Lemma 4.10]{fat1} this type is eventually transformed into the homaloidal type $(2;1^3,0^{r-3})$ after a finite iteration of arithmetic quadratic transformations.
By the same iteration, the vector $(s+1;2\mu_4,\ldots,2\mu_r,2^3)$ gets transformed into the vector $(2;2^3,0^{r-3})$ and by Lemma~\ref{transfquadpreserva} applied iteratively the corresponding linear systems in degree $s+1$ and $2$, respectively, have the same vector dimension.
But since $\{2\mu_4,\ldots,2\mu_r,2^3\}$ are the virtual multiplicities of $\widetilde{J}^{(2)}$ then $\dim_k (\widetilde{J}^{(2)})_{s+1}=\dim_k I_{2}$,
where $I=p_1^2\cap p_2^2\cap p_3^2=(x,y)^2\cap (x,z)^2\cap (y,z)^2$.
Clearly, $\dim_k I_{2}=0$.

\medskip

(c) Set $N:=(s+3/2)$. 
By Proposition~\ref{tilde-generation}, the ideal $\widetilde{J}$ is linearly presented. 
Referring to item (a), since $\{xf, yf, xg,zg\}$ is $k$-linearly independent we can assume that they form a subset of minimal generators of $\widetilde{J}$.
For the sake of dealing with the presentation matrix of $\widetilde{J}$ we will write  $xf, yf, xg,zg$  as the first four generators in this order.
With this proviso, we may assume that its $(N+1)\times N$ presentation matrix $\varphi$ has the a concatenation shape
$$\varphi=(\varphi_1\;|\;\varphi_2)$$
where $\phi_1$ is the transpose of the matrix
$$\left(\begin{array}{cccccccc}
y&-x&0&0&0&\cdots&0\\
0&0&z&-x&0&\cdots&0
\end{array}\right).$$

{\sc Claim 1:} no $k$-linear combination of the columns of $\phi_2$ is a vector in either one of the following column spaces of column size $N+1$:
\begin{equation}
\label{V1}
\left(\begin{array}{ccccc|ccccc|cc}
x&0& 0&\cdots & 0  &y&0& 0&\cdots  & 0&0 & 0     \\
0 &x&0& \cdots & 0  &0 &y& 0& \cdots   & 0&0 &0     \\
0&0&x& \cdots & 0  &0 &0& y & \cdots &0  & z& 0\\
0&0&0&\cdots &0&    0 & 0  &0& \cdots & 0  & 0 & z\\
\vdots &\vdots & \vdots & \vdots & \vdots &\vdots & \vdots  &\vdots &\vdots & \vdots &\vdots &\vdots\\
0 &0 & 0& \cdots  &x & 0 &0 & 0 & 0 & y&0 &0
\end{array}\right)
\end{equation}
and
\begin{equation}
\label{V2}
\left(\begin{array}{ccccc|ccccc|cc}
x&0& 0&\cdots & 0  &y&0& 0&\cdots  & 0&z & 0     \\
0 &x&0& \cdots & 0  &0 &y& 0& \cdots   & 0&0 &z     \\
0&0&x& \cdots & 0  &0 &0& y & \cdots &0  & 0& 0\\
0&0&0&\cdots &0&    0 & 0  &0& \cdots & 0  & 0 & 0\\
\vdots &\vdots & \vdots & \vdots & \vdots &\vdots & \vdots  &\vdots &\vdots & \vdots &\vdots &\vdots\\
0 &0 & 0& \cdots  &x & 0 &0 & 0 & 0 & y&0 &0
\end{array}\right)
\end{equation}
We argue with the first of these spaces, the argument for the second one being analogous.
 
Assuming to the contrary, let a vector $\mathbf{v}$ in the first space be a $k$-linear combination of the columns of $\phi_2$. One can assume $\mathbf{v}$ is actually a column of $\phi_2$.
But $xg$ is the determinant of $\phi$ by striking out the third row, namely:

$$xg=\det \left(\begin{array}{cc|c}
y&0&\\
-x&0&\\
0&-x&\\
0&0&\mathbf{v}^-\\
\vdots&\vdots&\\
0&0&\\
\end{array}\right),
$$
where $\mathbf{v}^-$ denotes the vector $\mathbf{v}$ with its third coordinate omitted.
Expanding by Laplace along the second column and canceling $x$, one sees that $g$ belong to the ideal $(x,y)^2$. In other words, $g$ has multiplicity $\geq 2$ at the point $p_1=(0:0:1)$. This contradicts the fact that $g$, being the equation of a principal curve, its multiplicity at $p_1$ is $1$.

\smallskip

Henceforth, we assume that $\phi$ has the proposed shape.
Let $\tt:=\{t_1,\ldots,t_{N+1}\}$ denote coordinates of the target space $\pp^N$, where $N=(s+3)/2$.
Consider the incidence ideal $I_1(\tt\cdot \phi)\subset R[\tt]$.
Since $\phi$ is a linear matrix with $\phi_1$ the initial submatrix, the Jacobian matrix of these biforms with respect to $x,y,z$ has the following shape

$$B=\left(\begin{array}{ccc}
-t_2&t_1&0\\
-t_4&0&t_3\\
\ell_{3,1}&\ell_{3,2}&\ell_{3,3}\\
\ell_{4,1}&\ell_{4,2}&\ell_{4,3}\\
\vdots&\vdots&\vdots\\
\ell_{N,1}&\ell_{N,2}&\ell_{N,3}
\end{array}\right),$$
where $\ell_{i,j}$ is a certain linear form in $\tt$.

\smallskip

{\sc Claim 2:} each of the following two sets of $1$-forms in $k[\tt]$ is $k$-linearly independent
$$\{t_3,t_4,\ell_{3,3},\ldots,\ell_{N,3}\}
\quad \mbox{and}\quad
\{t_1,t_2,\ell_{3,2},\ldots,\ell_{N,2}\}$$

This will follow from the specified shape of $\phi=(\phi_1|\phi_2)$.
Namely, for $1\leq i\leq N+1$ and $3\leq j\leq N$, the $(i,j)$-entry  of $\varphi_2$ is

\begin{equation}\label{entradasphi_2}
\frac{\partial\ell_{j,1}}{\partial t_i}x+\frac{\partial\ell_{j,2}}{\partial t_i}y+\frac{\partial\ell_{j,3}}{\partial t_i}z
\end{equation}
Say we have
	\begin{equation}
	\alpha_1 t_3+\alpha_2 t_4=\sum_{u=3}^{N}\beta_u \ell_{u3}
	\end{equation} 
for certain $\alpha_1,\alpha_2,\beta_3,\ldots,\beta_{N}\in k.$

Taking derivatives with respect to the variables $\tt$ yelds
\begin{equation}\label{derivando}
\sum_{u=3}^{N}\beta_u \frac{\partial{\ell_{u3}}}{\partial t_i}=\left\{\begin{array}{cccc}
\alpha_1,& \mbox{if}\;i=3\\
\alpha_2,&\mbox{if}\;i=4\\
0,&\mbox{otherwise}
\end{array}\right.
\end{equation}

Therefore, it suffices to show that $\beta_u=0$ for $3\leq u\leq N$.
If we form a $k$-linear combination of the columns of $\phi_2$ with the $\beta_i$'s as coefficients, then \eqref{entradasphi_2} e \eqref{derivando} yield a vector in the column space (\ref{V1}); this contradicts the statement of Claim 1.

The argument for the other set will be analogous, using Claim 1 as regards the column space (\ref{V2}) instead.

\medskip

We now come towards the inverse map of $\mathfrak{F}$.
First, note that the image of the map cannot be cut by a linear form because the map is defined by a $k$-basis of a linear system.
Then, by  \cite[Theorem 2.18]{AHA}, the $2$-minors of any one among the following matrices gives a representative of the inverse map:
\begin{equation}\label{matr_definem_inversa}
\left(\begin{array}{cccc}-t_2&t_1&0\\ -t_4&0&t_3\end{array}\right),\quad \left(\begin{array}{cccc}-t_2&t_1&0\\ \ell_{i,1}&\ell_{i,2}&\ell_{i3}\end{array}\right)\quad \mbox{ou}\quad \left(\begin{array}{cccc}-t_4&0&t_3\\ \ell_{i,1}&\ell_{i,2}&\ell_{i3}\end{array}\right)
\end{equation}
where $3\leq i\leq N,$ the reason being that the rank of any of these matrices modulo the defining equations of the image is $2$, as the minors are products of linear forms and the defining equations of the image generate a prime ideal with initial degree $\geq 2$ by the above observation.

Each of these representatives is spanned in degree $2$ and there are $2(N-2)+1=2N-3=s$ of them.
At this point it is not at all obvious that no two of these are proportional, thus avoiding repetition.
To deal with this question, we will instead look at the corresponding (source) inversion factors as described before (see Section~\ref{basic-term}) and prove an even stronger result.

Let $D_1,\ldots, D_s\in R$ denote these factors.
Note that $\deg(D_i)=2N-1=s+2$.
By Proposition~\ref{inversionfactor_is_symbolic} they belong to the second symbolic power of $\widetilde{J}$.
Thus, we will be through provided we prove the following

\smallskip

{\sc Claim 3:} $D_1,\ldots, D_s$ are $k$-linearly independent.

\smallskip

First one observes that letting $\mathbf{h}:=\{xf,yf,xg,zg,h_5,\ldots,h_{N+1}\}$ denote a $k$-basis of $\widetilde{J}_{N}$, an easy calculation shows that the inversion factors corresponding to the three blocks of matrices as above have the respective forms
$$xfg,\quad f\ell_{i,3}({\bf h}),\quad \mbox{and}\quad g\ell_{i2}({\bf h})$$
where  $3\leq i\leq N.$

Say we have
	$$\alpha xfg+f\sum_{i=3}^{N}\beta_i\ell_{i,3}({\bf h})+g\sum_{i=3}^{N}\gamma_i\ell_{i,2}({\bf h})=0,$$ 
	where $\alpha,\beta_i,\gamma_i\in k.$ 
	
	From this, there are linear forms  $L$ e $L'$ de $k[x,y,z]$ such that 
	
	$$Lf=\displaystyle\sum_{i=3}^{N}\gamma_i\ell_{i,2}({\bf h})
	\quad\mbox{and}\quad
	L'g=\displaystyle\sum_{i=3}^{N}\beta_i\ell_{i,3}({\bf h}).$$ 
	
	Therefore, $Lf,L'g\in \widetilde{J},$ and hence  $L=a_1x+b_1y$ and $L'=a'_{1}x+b'_{1}z$ for certain $a_1,a'_1,b_1,b'_1\in k.$ 
	It then follows that
	$$a_1t_1+b_1t_2-\sum_{i=3}^{N}\gamma_i\ell_{i,2}(\tt) \quad \mbox{and}\quad  a'_1t_3+b'_1t_4-\sum\beta_i\ell_{i,3}(\tt)$$
	are among the defining equations of the image of $\mathfrak{F}$.
	However, as remarked earlier, the image is not cut by any linear form.
	Consequently,
		$$a_1t_1+b_1t_2-\sum_{i=3}^{N}\gamma_i\ell_{i,2}(\tt)=a'_1t_3+b'_1t_4-\sum\beta_i\ell_{i,3}(\tt)=0.$$
	By Claim 2, $\beta_u=0$ and $\gamma_u=0$, as was to be shown.
	\qed
	
\begin{Remark}\rm
Since $\widetilde{J}$ is a linearly presented ideal, item (a) of the proposition would follow by \cite[Theorem 3.2]{AHA} provided we knew a priori that the $k$-subalgebra $k[\widetilde{J}_N]\subset R$ has dimension $3$ (i.e., that the map is dominant). 
\end{Remark}
  
  \medskip
  
  With the above notation, we file the following

\begin{Conjecture}\label{structure-second-symb} 
\begin{enumerate}
\item[{\rm (a)}] $\widetilde{J}^{(2)}=(D_1,\ldots,D_s).$ 
\item[{\rm (b)}] $\widetilde{J}^2=(x,y,z)\widetilde{J}^{(2)}=(xf,yf,xg,zg)\widetilde{J}.$ 
\end{enumerate}
\end{Conjecture}
\demo
Assuming item (a) valid, (b) can be shown in the following way,
First, the leftmost equality goes as follows.
By Corollary~\ref{tilde-symbolic-powers}~(b), one has $(\widetilde{J}^2)_{2N}=(\widetilde{J}^{(2)})_{2N}$.
By (a), $\widetilde{J}^{(2)}=(D_1,\ldots,D_s).$, hence $(\widetilde{J}^{(2)})_{2N-1}=\sum_i kD_i$.
It follows that $(\widetilde{J}^2)_{2N}$ is the span of $(x,y,z)_1\, (\widetilde{J}^{(2)})_{2N-1}$.
But since both $\widetilde{J}^2$ and $\widetilde{J}^{(2)}$ are generated in one single degree, we get the desired equality.

As for the rightmost equality, by definition $(x,y,z)D_i\, (1\leq i\leq 2N-1=s)$ is contained in the ideal of $2$-minors of one among the matrices in (\ref{matr_definem_inversa}), further evaluated on the generators $\mathbf{h}=\{xf,yf,xg,zg,h_5,\ldots,h_{N+1}\}$.
These evaluated matrices have the form
$$\left(\begin{array}{cccc}-yf&xf&0\\ -zg&0&xg\end{array}\right), \left(\begin{array}{cccc}-yf&xf&0\\ \ell_{i,1}({\bf h})&\ell_{i,2}({\bf h})&\ell_{i3}({\bf h})\end{array}\right), \left(\begin{array}{cccc}-zg&0&xg\\ \ell_{i,1}({\bf h})&\ell_{i,2}({\bf h})&\ell_{i3}({\bf h})\end{array}\right).
$$
Clearly, any $2$-minor above is contained in the product $(xf,yf,xg,zg)\widetilde{J}$.
Therefore, we obtain the following string of relations:
$$\widetilde{J}^2=(x,y,z)\widetilde{J}^{(2)}\subset (xf,yf,xg,zg)\widetilde{ J}\subset\widetilde{J}^2.$$
Thus, we are home.
\qed

\medskip

\subsection{Main theorem}

 As a piece of independent interest, in the footsteps of Question~\ref{main_queries}, we file the following result.
 It will be used in the proof of Theorem~\ref{first_main_theorem}.
 
 \begin{Proposition}\label{degree_of_image}
Let $J\subset R$ denote the fat ideal of general points with a sub-homaloidal multiplicity set in degree $s$. If $J=(J_s)$ then the image of the rational map defined by the linear system $J_s$ is a variety of degree $s/r$, where $r:=(k(J_s):k(R_s))$ is the field extension degree.
 \end{Proposition}
 \demo
 Let $(s;\boldsymbol{\mu})$ denote the given sub-homaloidal type and let $P_i\, (1\leq i\leq r)$ denote the respective defining prime ideals of the given general points.
 
 We will apply the formula of \cite[Theorem 6.6 (a)]{ram2}.
 
 For this matter, we first claim that  $e_{J_{P_i}}(R_{P_i})=\mu_i^2$, for every $i$ -- here, for an $\mathfrak{n}$-primary ideal $\mathfrak{a}$ in a local ring $(A,\mathfrak{n})$, the symbol $e_{\mathfrak{a}}(A)$ denotes the Samuel multiplicity of $A$ with respect to $\mathfrak{a}$.
 This follows by  recalling that $J_{P_i}={P_i}_{P_i}^{\mu_i}$ and applying Samuel's formula (see, e.g., \cite[Proposition 4.5.2 (b)]{BHbook}), with $\lambda (\_)$ denoting length; it obtains
 \begin{eqnarray}\label{Samuel}\nonumber
 e_{J_{P_i}}(R_{P_i})& = & \lim_{n\rar \infty}\frac{2}{n^2}\,\lambda\left(\frac{R_{P_i}}{{P_i}_{P_i}^{(n+1)\mu_i}}\right)=\lim_{n\rar \infty} \frac{2}{n^2}\,{{(n+1)\mu_i+1}\choose {2}}\\ 
 &=& \lim_{n\rar \infty} \frac{(n+1)^2\mu_i^2}{n^2}=\mu_i^2.
 \end{eqnarray}
 
 Now in the notation of \cite[Theorem 6.6 (a)]{ram2}  we take $R=k[x,y,z]$ to be our ground ring and $I=J$, an ideal of height $g=2$. Then $p$ runs through the primes $P_i$, i.e., the minimal primes of $R/J$.
 According to that result there is actually an equality
 $$e(k[J_s])= \frac{1}{r} \left(e(R)s^2-\sum_i \, e_{J_{P_i}}(R_{P_i})\, e(R/P_i)   \right)= \frac{1}{r} (s^2-\sum_i\mu_i^2),$$
 by (\ref{Samuel}) above and because $P_i$ is generated by linear forms.
 
 By the second equation of condition for $J$, one has $\sum_i\mu_i^2=s^2-s$.
 It follows that $e(k[J_s])= s/r$.
 \qed
 
 \medskip

 We will draw on the following result of Cortadellas--Zarzuela.
 For convenience, we state the part of their result in the form that best suits our present application:
 
 \begin{Theorem} {\rm (\cite[Theorem 5.12]{CorZar})}\label{CortZar}
 	Let $R$ denote a Cohen--Macaulay local ring and let $J\subset R$ stand for an ideal satisfying the following properties:
 	\begin{enumerate}
 		\item[{\rm (1)}]  $J$ is unmixed, and $R/J$ and $R/J^2$ have different depths.
 		\item[{\rm (2)}]  $\max\{r(J_P)\,|\, J\subset P, \hht(P)=\hht(J)\}\leq 1$.
 		\item[{\rm (3)}] $\ell(J)=\hht(J)+1\geq 3$.
 		\item[{\rm (4)}] $r(J)\leq 2$.
 	\end{enumerate}
 	Then ${\rm depth}\, \mathcal{R}_R(J)=\min\{{\rm depth} \,R/J, {\rm depth} \,R/J^2+1\} +\hht(J)+1$.
 \end{Theorem}
 Here $\ell(\_)$ denotes analytic spread, while $r(\_)$ stands for the minimum of all reduction numbers with respect to  minimal reductions.

Next is one of our main results about the sub-homaloidal virtual multiplicities introduced in the previous subsection.
It will give  in a special case an affirmative answer to all queries stated in Question~\ref{main_queries}.

It is easy to see that for a sub-homaloidal multiplicity set in degree $s$, the integer $s$ is necessarily odd.
Thus, we may assume that $s\geq 3$ as in the statement below.

\begin{Theorem}\label{first_main_theorem} 
Let $\boldsymbol\mu=\{\mu_1\geq\cdots \geq\mu_r\}$ denote a sub-homaloidal multiplicity set in degree $s\geq 3$.
Letting $\mathbf{p}\subset \pp^2$ denote a set of $r$ general points, write $J:=I(\mathbf{p}, \boldsymbol\mu)$ for the corresponding ideal of sub-homaloidal type.
Suppose that:
\begin{enumerate}
\item[{\rm (i)}]
The {\rm ``}doubled$\,${\rm ''} homaloidal type $(2s-1;2\boldsymbol\mu)$ is proper.
\item[{\rm (ii)}] $\mu_1=\mu_2=\mu_3= (s-1)/2$.
\end{enumerate}
Then:
\begin{enumerate}
\item[{\rm (a)}] $J$ is generated in degree $s$ and its minimal graded free resolution of $J$ has the form
{\large
	\begin{equation}\label{res_halved_mults}
	0\rar R^3\left(-\left(s+1\right)\right)\oplus R^{\frac{s-3}{2}}\left(-\left(s+2\right)\right)\stackrel{\varphi}\rar R^{\frac{s+5}{2}}\left(-s\right)\rar J\rar 0.
	\end{equation}}
\item[{\rm (b)}] The linear system $J_s$ defines a birational map of $\pp^2$ onto the image and the latter is a variety of degree $s$ in $\pp^{\frac{s+3}{2}}$.
\item[{\rm (c)}] The image of the birational map defined by $J_s$ is an arithmetically Cohen--Macaulay
surface in its  embedding in $\pp^{\frac{s+3}{2}}$ and, moreover, it is ideal theoretically defined by quadrics and cubics.
\end{enumerate}
\end{Theorem}

\demo  Since the points in $\mathbf{p}$ are general there is a Cremona map with homaloidal type $(2s-1;2\boldsymbol\mu)$ and these base points.

(a) That the initial degree of $J$ is $s$ is the content of Lemma~\ref{indeg_of_J} -- note that this result depended on the fact that
$\mathbb{J}:=I(\mathbf{p},2 \boldsymbol\mu)$ is the second symbolic power of $J$ and on the equations of sub-homaloidness, but not on assumption {\rm (}i{\rm )} about the three highest virtual multiplicities.

To see that $J$ is actually generated in degree $s$ under the $3$-uniform hypothesis of (ii) we  draw upon Corollary~\ref{equi-generation} and, moreover, from the proof of Proposition~\ref{tilde-generation} (a), we know that $\dim J_s= \frac{s+5}{2}$, the expected dimension.

Finally, $J$ is a codimension $2$ perfect ideal, the minimal free resolution of $J$ has the form $$0\rar \oplus_i R^{u_i}(-i)\rar
\oplus_i R^{v_i}(-i)=R^{g}(-s)\rar J\rar 0$$
with $g=(s+5)/2,$  e $i> s.$
The Betti numbers $\nu_i$ can be computed as  $$u_i=v_i-\Delta^{3}h_J(i),$$
where $\Delta$ denotes the difference function and $h_{J}$ the Hilbert function of $J$ (see, e.g., \cite{FiHaHo}).  Here then $$u_i=\nu_i-h_J(i)+3h_J(i-1)-3h_J(i-2)+h_J(i-3),$$
which with the  data coming from item (i) yields

$$u_{s+1}=-h_J(s+1)+3h_J(s)-3h_J(s-1)+h_J(s-2)=-3\frac{(s+3)}{2}+3
\frac{s+5}{2}=3$$
and
$$u_{s+2}=-h_J(s+2)+3h_J(s+1)-3h_J(s)+h_J(s-1)=-5\frac{s+3}{2}+9
\frac{s+3}{2}-3\frac{(s+5)}{2}=\frac{s-3}{2}.$$

Since $u_{s+1}+u_{s+2}=g-1,$ it obtains
$$\oplus R^{u_i}(-i)= R^{3}(-(s+1))\oplus R^{g-4}(-(s+2)),$$
as claimed.

\medskip

(b) By Lemma~\ref{transfquadpreserva}, the rational map defined by the linear system $J_s$ is obtained from the rational map defined by the linear system of $\widetilde{J}_{\widetilde{s}}$.
Therefore, the result follows from Theorem~\ref{birr_tildeJ} (a).
	
A alternative direct proof can be given in the same pattern as in the proof of Theorem~\ref{birr_tildeJ} (a).
For this purpose we consider the following ideals of fat points:
\begin{equation}\label{jota1}
J'=P_1^{(s-3)/2}\cap P_2^{(s-1)/2}\cap P_3^{(s-1)/2}\cap P_4^{\mu_4}\ldots\cap P_r^{\mu_r},
\end{equation}

\begin{equation}\label{jota2}
J''=P_1^{(s-1)/2}\cap P_2^{(s-3)/2}\cap P_3^{(s-1)/2}\cap P_4^{\mu_4}\ldots\cap P_r^{\mu_r}
\end{equation}

By \cite[Lemma 4.10]{fat1}, the corresponding types are proper exceptional types in degree $s-1$ in the sense of \cite[Definition 5.5.1]{alberich}.
Therefore, by \cite[Proposition 5.5.13]{alberich} there exist irreducible forms $f\in J'_{s-1}$ and $g\in J''_{s-1}$.
We can assume that $P_1=(x,y)$ and $P_2=(x,z)$.
Consider the linear subsystem $L\subset J_s$ spanned by the $s$-forms $xf,yf,xg,zg$.
Clearly then the ideal $(L)\subset R$ has codimension $2$.
Since $s\geq 3$, $L$ defines a birational image onto its image by Lemma~\ref{4guys}. Then the larger linear system $J_s$ also defines a birational map of $\pp^2$ onto the image.

This proves the birationality part of the item.

The second statement of the item follows from Proposition~\ref{degree_of_image} or, alternatively, from Proposition~\ref{image_has_degree_s}.

\medskip

(c) We first deal with the Cohen--Macaulay assertion, which means that the homogeneous coordinate ring of the image is a Cohen--Macaulay ring.
For this, note that the latter is identified with the $k$-subalgebra
$k[J_s]\subset R$ up to degree rescaling.
The trick is to observe that $k[J_s]\simeq k[\tilde{J}_N]$ as graded $k$-algebras, where $N=(s+3)/2$, as an easy consequence of the material in Subsection~\ref{AQT}.

We apply \cite[Theorem 2.4]{GGP} taking $I$ in this theorem to be our $\tilde{J}$.
Observe that for this, one has to make sure that the scheme $Z$ corresponding to $\tilde{J}$ satisfies the property that no line of $\pp^2$ intersects $Z$ in more than ${\rm reg}(\tilde{J})$ points degree-wise, as is required in \cite[Theorem 2.1 (ii)]{GGP}.
But since we are taking general points, this degree count is at most $\mu_4+\mu_5$.
On the other hand, $(2s-1;s-1,s-1,s-1,2\mu_4,2\mu_5,\ldots, 2\mu_r)$ is a proper homaloidal type, hence so is its arithmetic transform $(s+1;2\mu_4,2\mu_5,\ldots, 2\mu_r, 1^3)$.
But for a proper homaloidal type one must have $2\mu_4+2\mu_5\leq s+1$, and so $\mu_4+\mu_5\leq (s+1)/2<(s+3)/2={\rm reg}(\tilde{J})$.

Finally, the corresponding divisor ($D_{\sigma}$ in the notation there, where $\sigma={\rm reg}(\tilde{J})$) is then very ample. But the image of the embedding it defines is the same as the image of the birational map defined by the linear system $\tilde{J}_N$ (see Theorem~\ref{birr_tildeJ} (a) for this map).
Therefore, we are through.

\smallskip

In order to argue about the generating degrees of the defining ideal of the image we will prove that the Rees algebra $\mathcal{R}_R(J)=R[Jt]\subset R[t]$ is Cohen--Macaulay.
Once this granted, we apply \cite[Corollary 6.3]{AHT}) by which the relation type of $J$ is at most $\dim R=3$ (the relation type is the highest presentation degree of a minimal homogeneous generator of a presentation ideal of $\mathcal{R}_R(J)$.)
But since the defining ideal of the image in the embedding is contained in $\mathcal{J}$, its minimal generators must be minimal generators of the latter. This implies that the image is ideal theoretically generated in degree $\leq 3$, hence is generated by quadrics and cubics.

To show the Cohen--Macaulayness of $\mathcal{R}_R(J)$ we will apply Theorem~\ref{CortZar}.
In order to apply it to our present ideal $J$, we take $R=k[x,y,z]_{(x,y,z)}$ as usual and go through the list of hypotheses of the theorem.
Condition (1)  has two parts. The first is obvious since $R/J$ is even Cohen--Macaulay; in particular, $R/J$ has depth $1$. The second statement follows from the fact that the symbolic power $J^{(2)}=\mathbb{J}\neq J^2$ as it contains elements of degree $2s-1$, while $J^2$ is generated in degree $\geq 2s$. Therefore, $R/J^2$ has depth zero.

Condition (2) goes as follows: locally at a minimal prime $P$ of $R/J$ the ideal $J$ is a power of certain order $\mu$ of a complete intersection which, up to a change of coordinates, can be thought of as $(x,y)^{\mu}$ in $k[x,y]_(x,y)$.
Clearly, $(x^{\mu},y^{\mu})$ is a minimal reduction of the power with reduction number $1$.
Thus the required maximum is attained at $1$.

\smallskip

Condition (3) is equivalent to having $\dim k[J_s]=3$, which follows immediately from the birationality assertion of item (b).

\smallskip

We now deal with condition (4), which is by far the hardest.
It will be a consequence of the Cohen--Macaulayness of the image via the following general result:

\begin{Proposition}\label{red_at_most_2} Suppose that $I\subset R=k[x,y,z]$ is a homogeneous ideal such that its special fiber $\mathcal{F}(I)$ is Cohen-Macaulay and has dimension $3$. The following conditions are equivalent:
	\begin{enumerate}
		\item[{\rm (i)}] The reduction number $r(I)$ is at most $2$.
		\item[{\rm (ii)}] The Hilbert function $H_{\mathcal{F}(I)}(t)$ of $\mathcal{F}(I)$ coincides with its Hilbert polynomial for all $t\geq 0$.
	\end{enumerate}
\end{Proposition}
\demo
Since the degree of the numerator of the Hilbert series of $\mathcal{F}(I)$ is at most $2$ by \cite[Proposition 1.85]{Wolmbook3}, the equivalence of (i) and (iii) follows by \cite[Proposition 4.1.12]{BHbook} or \cite[Proposition 4.3.5 (c)]{BHbook}.
\qed

\smallskip

The required estimate for the reduction number now follows from Proposition~\ref{image_has_degree_s}.

\medskip

Now apply the theorem with the present values, thus getting:
$${\rm depth} \mathcal{R}_R(J)=\min\{ 1, 0+1\} +2+1=4.$$
Since $\dim \mathcal{R}_R(J)=\dim R+1=4$, we have proved that the Rees algebra of $J$ is Cohen--Macaulay as claimed.
\qed

\begin{Remark}\rm
(1) One notes that \cite[Theorem 2.1 (ii)]{GGP} is not applicable above as was in \cite[Remark 2.7]{fat1} since, from (\ref{res_halved_mults}) the regularity of $J$ is $s+2-1=s+1$. Of course, once one has that $J_s$ defines a birational map of $\pp^2$ onto the image, then so does $J_{s+1}$.
Further, since $(2s-1;2\boldsymbol\mu)$ is assumed to be proper, one has $2(\mu_1+\mu_2)=2\mu_1+2\mu_2\leq 2s-1$, i.e., $\mu_1+\mu_2<s\, (<s+1)$; therefore, the extra condition in \cite[Theorem 2.1 (ii)]{GGP} is fulfilled, and hence the image of the birational map defined by $J_{s+1}$ is in fact smooth and arithmetically Cohen--Macaulay.
The question remains as to how, if at all,  the two images relate to each other -- a gloomy perspective as the image under $J_s$ has cubic minimal defining equations, while the one under $J_{s+1}$ seems to be defined by quadrics.

(2) Neither the statement of the above theorem nor its proof rests on a  hocus-pocus deformation to some ``general'' setup.
That is, if one lets $J\subset R$ denote a codimension $2$ perfect ideal with resolution as in (\ref{res_halved_mults}) and general matrix entries then the reduction number of $J$ can be significantly higher than the expected value $2=\ell(J)-\hht(J)+1$. This possibly explained by the fact that, although the map  is still birational onto its image, the latter is a non arithmetically Cohen--Macaulay image.
\end{Remark}

We close this subsection with the following result which seems to have gone unnoticed (see \cite[Proposition]{verde} for a more encompassing statement).

\begin{Proposition}\label{ReesCM-vs-fiberCM}
	Let $I\subset R=k[x,y,z]$ denote a homogeneous ideal generated in a single degree $s$ such that $\dim k[I_s]=3$. If the Rees algebra $R[It]\subset R[t]$ of $I$ is Cohen--Macaulay then so is the algebra $k[I_s]$.
\end{Proposition}
\demo Since the algebra $k[I_s]$ is graded isomorphic to the $k$-subalgebra $k[It]\subset R[It]$ it suffices to deal with the latter.
Since $I=(I_s)$ then $k[It]$ is a direct summand of $R[It]$ as $k[It]$-module. Let $Q\subset I$ denote a minimal reduction of $I$.
Such reductions exist in this setup that are generated by forms of the same degree as the degree of the linear system $I_s$ (an easy application of Noether normalization lemma).
Thus, we assume this is the case of $Q$.
We know that $Q$ is generated by analytically independent elements $q_1,q_2,q_3$ (with respect to the irrelevant ideal $\fm\subset R$), hence $k[Qt]\subset R[Qt]$ is a polynomial ring in $q_1t,q_2t,q_3t$ over $k$.
Set $h:=(q_1t,q_2t,q_3t)\subset k[Qt]$, an ideal of codimension $3$ in both $k[Qt]$ and $R[Qt]$.
Now, one has the two splittings $R[Qt]=k[Qt]\oplus \fm R[Qt]$ and $R[It]=k[It]\oplus \fm R[It]$ as $k[Qt]$-modules and as $k[It]$-modules, respectively.
In addition, one has a naturally induced integral ring extension $R[Qt]\subset R[It]$.
Thus, the image of $h$ under this extension also has codimension $3$ in $R[It]$.
Since the latter is assumed to be Cohen--Macaulay, the images of the elements $q_1t,q_2t,q_3t$  form a regular sequence in $R[It]$.
But, these elements belong to the direct summand $k[It]$, hence they form a regular sequence in the latter as well.
On the other hand, by assumption $\dim k[It]=\dim k[J_s]=3$.
Therefore, $k[J_s]$ is Cohen--Macaulay as required.
\qed

\begin{Remark}\rm
If one had a proof of the bound $r(J)\leq 2$ running independently of knowing that the special fiber is Cohen--Macaulay then one could apply the Cortadellas--Zarzuela criterion along with the above result to actually deduce the Cohen--Macaulayness of the special fiber.
\end{Remark}

\subsection{Structure of the Cremona base ideal}

In this part one relates the previous results back to the structure of the base ideal of the original Cremona map that motivated them, aiming
ultimately at the minimal free resolution of the base ideal.
Although the shape of the resolution had been dealt with in \cite[Propostion 4.12]{fat1}, the method employed here tells how birational invariants come into the picture.

\begin{Theorem}\label{second_main_theorem}
Keeping the notation and the hypotheses in the statement of {\rm Theorem~\ref{first_main_theorem}}, let $\mathfrak{F}:\pp^2\dasharrow \pp^{(s+3)/2}$ denote the birational map of $\pp^2$ onto the image defined by the linear system $J_s$.
Set $\pp^{(s+3)/2}={\rm Proj}(k[\tt])$, where $\tt=\{t_1,\ldots,t_{(s+5)/2}\}$.
\begin{enumerate}
\item[{\rm (a)}] Letting $\phi_1$ denote the $((s+5)/2)\times 3$ linear submatrix of $\phi$ as in {\rm (\ref{res_halved_mults})}, let $\rho$ stand for the $3\times 3$ Jacobian matrix of the three entries of $\tt\cdot \phi_1$ with respect to the ground variables $x,y,z$.
Then the inverse map $\mathfrak{F}^{-1}$ is defined by the ordered signed $2$-minors of any  $2\times 3$ submatrix of $\rho$ taken modulo the homogeneous defining ideal of the image of $\mathfrak{F}$.

\item[{\rm (b)}] The base ideal $I:=(\mathbb{J}_{2s-1})$ of the Cremona map with the associated homaloidal type $(2s-1; 2\boldsymbol\mu)$ is generated by the source inversion factors corresponding each to one of the three representatives of $\mathfrak{F}^{-1}$, as established in {\rm (}a{\rm )} and has a graded minimal pure resolution of the form
{\small\begin{equation}
0\rar R(-3s) \rar R^3(-(3s-1))\rar R^3(-(2s-1))\rar I\rar  0.
\end{equation}}
\end{enumerate}
\end{Theorem}
\demo
(a) The proof is sufficiently analogous to that of Proposition~\ref{structure-second-symb}~(c).
For this we use the construction in the alternative proof of Proposition~\ref{first_main_theorem}~(b), where one obtained $k$-linearly independent forms $\{xf,yf,xg,zg\}\subset J_s$. One can take these forms to be part of a basis of $J_s$, in this order; say, $\hh:=\{xf,yf,xg,zg, h_5,\ldots, h_{(s+5)/2}\}$.
Thus, we may assume that $\phi=(\phi_1|\phi_2)$ where $\phi_1$ is the transpose of the $3\times s$ matrix
$$\left(\begin{array}{cccccccc}
y&-x&0&0&0&\cdots&0\\
0&0&z&-x&0&\cdots&0\\
\ell_1&\ell_2&\ell_3&\ell_4 &\ell_5 &\ldots &\ell_s
\end{array}\right),$$
for certain linear forms $\ell_i, 1\leq i\leq s$, two of which at least are nonzero.

Thus, $\rho$ has the following shape

$$B=\left(\begin{array}{ccc}
-t_2&t_1&0\\
-t_4&0&t_3\\
L_1(\tt)&L_2(\tt)&L_3(\tt)\\
\end{array}\right),$$
where $L_i(\tt), i=1,2,3,$ are certain linear forms in $\tt$ two of which are non-vanishing.

{\sc Claim:} up to elementary column operations on $\phi_1$, one can assume that $L_1(\tt)=0$.

Note that the assertion is equivalent of showing that up to elementary column operations on $\phi_1$ every $\ell_i\in (y,z)$.

We argue as follows.
We first claim indeed $\ell_i\in (y,z)$ for any $i\neq 2,4$. 
Namely, for such a value of $i$ the corresponding $\ell_i$ appears on one of the following $3\times 3$ matrices
$$\left(\begin{array}{ccc}
y&0&\ell_1\\
-x&0&\ell_2\\
0&-x&\ell_4
\end{array}\right)\quad \mbox{if}\quad i=1,$$

$$\left(\begin{array}{ccc}
-x&0&\ell_2\\
0&z&\ell_3\\
0&-x&\ell_4
\end{array}\right)\quad \mbox{if}\quad i=3$$ ou

$$\left(\begin{array}{ccc}
-x&0&\ell_2\\
0&-x&\ell_4\\
0&0&\ell_{i}
\end{array}\right)\quad \mbox{if}\quad i>4.$$
Let $\Delta$ denote the determinant of any of these matrices.
It is easy to verify that $\Delta\in (y,z)$ if and only if $\ell_i\in (y,z)$.
Thus, assuming $\ell_i\notin (y,z)$, one can invert $\Delta$ in the localization $k[x,y,z]_{(y,z)},$ and hence the ideal $J_{(y,z)}$ is the ideal of maximal minors of a Hilbert-Burch matrix with $(s+3)/2-3=(s-3)/2$ columns.
This means that $J_{(y,z)}$  is generated by $(s-3)/2+1=(s-1)/2$ elements.
But this contradicts the given data where
$$J_{(y,z)}=(y,z)_{(y,z)}^{(s-1)/2}$$
is minimally generated by $(s-1)/2+1$ elements. 

Finally, as for the $x$-term in both $\ell_2$ and $\ell_4$ it can be canceled by obvious column operations using the other two linear syzygies

As a result, the earlier matrix has now the form
$$B=\left(\begin{array}{ccc}
-t_2&t_1&0\\
-t_4&0&t_3\\
0&L_2(\tt)&L_3(\tt)\\
\end{array}\right),$$
where $L_2(\tt)$ and $L_3(\tt)$ are non-vanishing linear forms in $\tt$.

Therefore, each of the $2\times 3$ submatrices thereon have a non-vanishing $2$-minor modulo the homogeneous ideal of the image of $\mathfrak{F}$, the reason being that these minors are all products of linear forms while the latter ideal is prime and is not cut by linear forms (since $\mathfrak{F}$ is defined by $k$-linearly independent forms). 

\smallskip

(b) Let $D_1,D_2,D_3$ denote the source inversion factors corresponding, respectively, to the representatives of  $\mathfrak{F}^{-1}$ given by the $2$-minors of the matrices
$$\left(\begin{array}{ccc}
-t_2&t_1&0\\
-t_4&0&t_3\\
\end{array}\right), \quad
\left(\begin{array}{ccc}
-t_2&t_1&0\\
0&L_2(\tt)&L_3(\tt)\\
\end{array}\right)\quad\mbox{and}
\quad
\left(\begin{array}{ccc}
-t_4&0&t_3\\
0&L_2(\tt)&L_3(\tt)\\
\end{array}\right)$$

By the definition of a source inversion factor, evaluating $\tt\mapsto \hh$ on the leftmost matrix gives
$$yD_1=\det\left(\begin{array}{ccc}
-yf&xf&0\\
-zg&0&xg\\
\end{array}\right)=-xyfg,
$$
hence $D_1=xfg$.

Proceeding similarly with the central matrix yields
$$yD_2=\det\left(\begin{array}{cc}
	-yf&0\\
	0&L_3(\hh)\\
\end{array}\right)=-yfL_3(\hh),
$$
hence $D_2=fL_3(\hh)$.
Calculating $D_2$ once more:
$$zD_2=\det\left(\begin{array}{cc}
-yf&xf\\
0&L_2(\hh)\\
\end{array}\right)=-yfL_2(\hh).
$$
Substituting yields $-yL_2(\hh)=zL_3(\hh)$.
Therefore, there exists a form $q\in k[x,y,z]$ of degree  $s-1$ such that $L_{3}(\hh)=yq$ (and forcefully, $L_2(\hh)=-zq$).

Thus, $D_2=yfq$.

Finally, calculating $D_3$ from the rightmost matrix gives
$$yD_3=\det\left(\begin{array}{cc}
-zg&xg\\
0&L_3(\hh)
\end{array}\right)=-zgL_3(\hh).$$
Substituting for $L_3(\hh)$ from above yields $D_3=zgq$.

Now that we have calculated the three source inversion factors, we claim that the form $q$ is the equation of the unique principal curve of degree $s-1$ with virtual multiplicities $((s-1)/2, (s-1)/2, (s-3)/2, \mu_4,\ldots,\mu_r)$.

To see this, first note that as above $zq=-L_2(\hh)\in J.$ 
Therefore,   $q$ and $zq$ will have the same multiplicity at any point off the line $\{z=0\}.$ 
But since the points are general  $q$ and $zq$ will have the same multiplicity at all points except possibly the coordinate points $p_2=\{x=z=0\}$ e $p_3=\{y=z=0\}.$ In addition, since $zq\in J$ the multiplicity of $q$ at $p_3$ is at least $(s-3)/2.$ 
Thus, to close the argument we need to show that the virtual multiplicity of $q$ at $p_2$ is   $(s-1)/2.$
But this follows from having $yq=L_3(\hh)\in J$ and the fact that $p_2$ does not lie on the line $\{y=0\}.$

\smallskip

By (a), the inverse of $\mathfrak{F}$ is defined by $2$-forms.
Then, by definition, each $D_i$ has degree $2s-1$.
By Proposition~\ref{inversionfactor_is_symbolic}, $D_i\in J^{(2)}$, for $i=1,2,3$. Since $J^{(2)}_{2s-1}$ generates the base ideal of a Cremona map of degree $2s-1$ it suffices to show that $\{D_1,D_2,D_3\}$ is $k$-linearly independent.
But drawing on the above explicit shape of these generators and using that $f,g,q$ are irreducible one easily verifies the assertion.

The proof of the statement concerning the  minimal free resolution of $I$ is the same as in \cite[4.12]{fat1}.
\qed

\begin{Remark}\rm
Note that for such a minimal $(2s-1)$-linear resolution as in (iii) above, $3s-1$ is the smallest possible twist. This can be verified in different ways -- see, e.g.,  \cite[Section 1]{HS} for a detailed account involving bounds for the regularity of $R/J$, or \cite[Corollary 1.2]{ST} for an elementary argument on degrees of syzygies.
\end{Remark}

We end this part by observing that the method of proof of Proposition~\ref{second_main_theorem}~(a) leaves a bit to desire in the sense that it involves the third principal curve $q$ a posteriori, while it would seem more balanced to start {\em ab initio} with the three $f,g,q$.
If one does not care about the double nature of these forms as equations of principal curves and as irreducible components of source inversion factors of a birational map, then one can proceed as in the proof of \cite[4.12]{fat1}.

Dealing from the beginning with the three events would lead us to consider the $6$ forms $\{xf,yf,xg,zg,yq,zq\}\subset J$ and retrace all arguments {\em ab initio}, such as the proof of Proposition~\ref{tilde-generation}~(c) and Lemma~\ref{tiny-res}.
As to Lemma~\ref{tiny-res} the question arises as to whether these forms are always linearly independent over $k$ (equivalently, whether $f,g,q$ do not admit a linear syzygy). Luckily, one can still conclude that the map defined by these forms is birational even in the case where the image is cut by a linear form.
However, reworking the proof of Proposition~\ref{tilde-generation}~(c) -- and, similarly, the one of Proposition~\ref{second_main_theorem}~(a) -- would face some technical difficulty with the linear dependence of the six forms.

In order to have a clean formulation and neat proofs, it would look like the following conjectural question ought to be settled beforehand.
The ground hypotheses are that the subhomaloidal type is $(s;(s-1)/2^3,\mu_4,\ldots,\mu_r)$ with corresponding proper homaloidal type, and $f,g,q$ denote the  principal curves as described above.

\begin{Conjecture}\rm
The ideal $(f,g,q)\subset R=k[x,y,z]$ is perfect (i.e., unmixed).
In addition, this ideal admits a linear syzygy if and only if either $s=3,5$ or else $s>5$ and $\mu_4=(s-1)/2$.
\end{Conjecture} 
We believe that the perfectness of the ideal $(f,g,q)$ is an algebraic expression of the geometric symmetry behind the class under consideration.
As for the exceptional subclass in the statement, see Example~\ref{special-type}.

\subsection{Relation to Bordiga--White parametrizations}

The notion of a Bordiga--White surface has been introduced in \cite[Definition 3.1]{GGP}.
We wish to relate this class of surfaces to the previous results.

As in the previous section, let $\boldsymbol\mu=(\mu_1,\ldots,\mu_r)$ denote a sub-homaloidal multiplicity set in degree $s\geq 3$ ($s$ necessarily odd)  and such that the associated doubled homaloidal type $(2s-1;2\boldsymbol\mu)$ is proper.
Letting $\mathbf{p}\subset \pp^2$ denote a set of $r$ general points, fix a plane Cremona map $\mathfrak{F}$ with this type and further write $J:=I(\mathbf{p}, \boldsymbol\mu)$.

Let us assume that $J$ is generated in degree $s$ (cf. Question~\ref{main_queries}, (b)) -- as is the case, e.g., under the hypotheses of Theorem~\ref{first_main_theorem}.
Let as well $\widetilde{J}\subset R$ as in in Subsection~\ref{AQT} denote the ideal of the transformed points $\{p_1,p_2,p_3,\mathcal{Q}(p_4),\ldots, \mathcal{Q}(p_r)\}$ with the virtual multiplicities of $S(\mathcal{Q})$

A general partial result away from the strict hypothesis (ii) of Theorem~\ref{first_main_theorem} goes as follows.

\begin{Proposition}\label{Bordiga_main}
Assume that the linear system $J_s$ has the expected dimension and that $\widetilde{J}$ is generated in degree $d:=2s-(\mu_1+\mu_2+\mu_3)$, where  $\mu_1\geq\mu_2\geq\mu_3$ are the three highest multiplicities. The following conditions are equivalent:
\begin{enumerate}
	\item[{\rm (i)}] $\widetilde{J}$ is minimally generated by $d+1$ forms.
	\item[{\rm (ii)}] $e(R/\widetilde{J})={{d+1}\choose {2}}$.
	\item[{\rm (iii)}]  $\mu_1+\mu_2+\mu_3 = 3 (s-1)/2$.
\end{enumerate}
Moreover, under any of these equivalent conditions $\widetilde{J}_d$ defines a birational map to $\pp^d$ whose image  is a Bordiga--White surface.
\end{Proposition}
\demo
(i) $\Rightarrow$ (ii)
Since $\widetilde{J}$ is generated in degree $d$ and $R/\widetilde{J}$ is Cohen--Macaulay, it is generated by the maximal minors of a $(d+1)\times d$ matrix with homogeneous entries. It follows that all entries must be linear forms.
This implies that $\widetilde{J}$ has minimal free resolution of the form
\begin{equation}\label{tiny-res}
0\rar R^{d}(-(d+1))\lar R^{d+1}(-d)\rar \widetilde{J}\rar 0.
\end{equation}
Computing the degree $e(R/\widetilde{J})$ out of the Hilbert series read on (\ref{tiny-res}), one easily gets the stated value in (ii).

\smallskip

(ii) $\Rightarrow$ (iii)
We express the degree $e(R/\widetilde{J})$  from its known in terms of the virtual multiplicities. 
By the facet of $\widetilde{J}$ as $J^{\mathcal{Q}}$, one has
$$m_1=s-(\mu_2+\mu_3),\; m_2=s-(\mu_1+\mu_3), \; m_3=s-(\mu_1+\mu_2),$$
yielding $m_1+m_2+m_3=3s-2\underline{\mu}=(2s-\underline{\mu})+ (s-\underline{\mu})= d+(s-\underline{\mu})$, where we have written $\underline{\mu}:=\mu_1+\mu_2+\mu_3$ for lighter reading.

Since $m_i=\mu_i$ for $i\geq 4$, it obtains:
\begin{eqnarray}\nonumber
\sum_{i\ge 1} m_i &=& m_1+m_2+m_3+\sum_{i\ge 4} m_i=d+(s-\underline{\mu})+ \sum_{i\ge 4}\mu_i\\ \nonumber
&=&d+(s-\underline{\mu})+ 3(s-1)- \underline{\mu}=d+2(2s-\underline{\mu})-3\\ \nonumber
&=& d+2d-3=3(d-1).
\end{eqnarray}
The computation of the sum os squares is slightly more involved but follows the same pattern.
As $m_1^2= s^2-2s(\mu_2+\mu_3)+\mu_2^2+\mu_3^2+2\mu_2\mu_3$, and $m_2^2,m_3^3$ have similar expressions, we derive
{\small
	\begin{eqnarray}\nonumber
	\sum_{i\ge 1} m_i^2 &=& m_1^2+m_2^2+m_3^2+\sum_{i\ge 4} \mu_i^2= m_1^2+m_2^2+m_3^2+s(s-1)-(\mu_1^2+\mu_2^2+\mu_3^2)\\ \nonumber
	&=& 3s^2-4s\underline{\mu}+2(\mu_1\mu_2+\mu_1\mu_3+\mu_2\mu_3) + 2(\mu_1^2+\mu_2^2+\mu_3^2)+s^2-s -(\mu_1^2+\mu_2^2+\mu_3^2) =\\ \nonumber
	&=& 4s^2-s(1+4\underline{\mu})+2(\mu_1\mu_2+\mu_1\mu_3+\mu_2\mu_3)+\mu_1^2+\mu_2^2+\mu_3^2 \\ \nonumber
	&=& 4s^2-4s\underline{\mu}+\underline{\mu}^2 -s=(2s-\underline{\mu})^2-s=d^2-s=d^2-(d+\underline{\mu})/2.
	\end{eqnarray}
}
Assembling yields
$$\frac{1}{2}(m_i(m_i+1))=\frac{1}{2}\left(3(d-1)+ d^2-\frac{d+\underline{\mu}}{2}\right)=\frac{1}{4}\left(d(2d+5)-(\underline{\mu}+6)\right).$$
By assumption, this value is also ${{d+1}\choose {2}}$.
This forces, upon an easy calculation, that $\underline{\mu}=3(d-2)$. Since $d=s-\underline{\mu}$, this gives in turn $\underline{\mu}=(3/2)(s-1)$.

\smallskip

(iii) $\Rightarrow$ (i)
 The expected dimension of $J_s$ is $\dim_k(J_s)=(s+5)/2=(s+3)+1$, while $d=2s-(\mu_1+\mu_2+\mu_3)= 2s-3(s-1)/2=(s+3)/2$.
 Since $\widetilde{J}$ is generated in degree $d=(s+3)/2$ and by construction we must have $\dim_k(J_s)=\dim_k(\widetilde{J}_{(s+3)/2})$ then the minimal number of generators of $\widetilde{J}$ is $d+1$.
 
\medskip

The birationality of the map follows from the criterion of \cite[Theorem 3.2]{AHA}, provided one shows that the dimension of the $k$-algebra $k[\widetilde{J}_d]$ is $3$ (maximal possible).
For this, recall once more the nature of $\widetilde{J}$ 
which came from the linear system  $J_s$ evaluated at the three algebraically independent $2$-forms defining $\mathcal{Q}$, further canceling a common factor of the resulting forms (Lemma~\ref{transfquadpreserva}).
Since $\dim k[J_s]=3$, clearly then $\dim k[\widetilde{J}_d]=3$ as well.

In order to conclude that the image is a Bordiga--White surface, we check the present situation against the ingredients of  \cite[Definition 3.1]{GGP}, whereby the ideal $I$ in this reference is taken to be $\widetilde{J}$.
First, from (\ref{tiny-res}) we see that the regularity of $R/\widetilde{J}$ is $d$.
Second, the degree of $R/\widetilde{J}$ is the assumed combinatorial number.
Finally, the condition on how a line intersects the corresponding ideal of fat points translates into the inequality $m_1+m_2<d$. The latter holds since $m_1+m_2=2s-(\mu_1+\mu_2+\mu_3)-\mu_3=d-\mu_3<d$ as $\mu_3>0$ due to Noether's inequality (\cite[Proposition 2.6.4]{alberich}).
\qed

\begin{Corollary}\label{main_prone_bis}
	Fix the notation and hypotheses of {\rm Theorem~\ref{first_main_theorem}}. Set $d:=(s+3)/2$ and let $\widetilde{J}$ be as in the previous theorem.
	Then statements {\rm (i)} through {\rm (ii)} and the birationality statement of the previous theorem hold true and, moreover, the image of the Bordiga--White map defined by  $\widetilde{J}_d$ has degree $2d-3$.
\end{Corollary}
\demo The assertions are clear since the hypotheses of the previous theorem are implied by Theorem~\ref{first_main_theorem}. To see the degree of the image, note that composing  with a quadratic transformation does not change the properties of the map.
Therefore, the image under $\widetilde{J}_d$ has the same degree as the image under $J_s$, which is $s=2d-3$ by Theorem~\ref{first_main_theorem} (b).
\qed

\begin{Remark}
	\rm
	The map defined by $\widetilde{J}_d$ lifts to a well-known embedding of $\pp^2$ blown-up on the base points of $\widetilde{J}_d$, induced by the corresponding very ample divisor, and the embedded surface is one of the Bordiga--White surfaces, studied in \cite{GGP}, where it is shown in addition that its homogeneous defining ideal of is generated by th $2$-minors of a suitable $2\times (d-1)$ matrix.
\end{Remark}

\section{Appendix: selected examples}

In this appendix, we give sufficient details of some examples illustrating parts of the theory.

\begin{Example}\rm
	The multiplicity set $(4,2^6, 1^2)$ is sub-homaloidal in degree $7$,  but the doubled homaloidal type $(13;8,4^6, 2^2)$ is not proper since it is eventually arithmetically transformed into the type $(4;-1,2^2,1^6)$.
	In the notation of Theorem~\ref{Bordiga_main}, $s=7$, $\underline{\mu}=8$ and $d=2\times 7-8=6$, while $(3/2)(s-1)=9>8$.
	The map defined by $\widetilde{J}_6$ is still birational with image a smooth arithmetically Cohen--Macaulay, but not a Bordiga--White surface for those invariants.
\end{Example}

\begin{Example}\rm
	The statement on birationality in Theorem~\ref{first_main_theorem} is still valid by relaxing condition (ii) to assume only that $\mu_1=\mu_2=(s-1)/2$, provided one can show that there exist forms $f\in J'_{s-1}$, $g\in J''_{s-1}$ such that $\gcd(f,g)=1$. 
	Unfortunately, in general this provision fails: in the sub-homaloidal type $(27;13^2,10,7^5,3^2,1)$ both $J'_{26}$ and $J''_{26}$ have dimension $1$ and any two respective nonzero equations have $\gcd$ of degree $\geq 4$ (this example has been obtained screening the output of a routine by A. Doria in {\em Singular}).
	One can also check that Proposition~\ref{tilde-generation} largely fails for the transformed ideal $\widetilde{J}$ associated to this example.
\end{Example}

\begin{Example}\label{special-type}\rm
Let $s\geq 3$ be an odd integer. For a given integer $r\geq 5$, consider the vector type $(2s-1; (s-1)^4,m_5,\ldots,m_r)$, where $s-1\geq m_5\geq\cdots\geq m_r$.
An easy calculation shows that this vector satisfies the equations of condition if and only if $m_5=\cdots=m_r=2$.
Thus, the given vector is homaloidal if and only if it is the vector $(2s-1; (s-1)^4,2^{s-1})$.
Moreover, the latter type is proper because applying an arithmetic quadratic transformation yields the homaloidal type $(s+1;s-1,2^{s-1},1^3)$ which is proper (see, e.g., \cite[Proposition 3.1]{fat1}). 
Now, the corresponding subhomaloidal type is $(s; ((s-1)/2)^4, 1^{s-1})$. Applying the same arithmetic quadratic transformation to the latter type gives the type $(d; d-2, 1^{2d-1})$, where $d:=(s+3)/2$. Letting $\widetilde{J}\subset R$ denote the ideal of fat general points with virtual multiplicities of this last type, then as an application of Corollary~\ref{main_prone_bis}, $\widetilde{J}$ is linearly presented, with  free resolution of the form
$$0\rar R^{d}(-(d+1))\lar R^{d+1}(-d)\rar \widetilde{J}\rar 0.$$
\end{Example}
We think that this is a very strict class of both homaloidal types and cases of Bordiga--White parameterizations, while the case where $m_1=m_2=m_3=s-1>m_4$ ought to be the ``generic'' component.
As evidence for this perception, we can also derive this example as an application of the following theorem:

\medskip

{\bf Theorem.} (\cite[Theorem A]{Dumnicki5}) {\em 
Let $\mathbf{p}=\{p_1,\ldots,p_r\}\subset \pp^2$ be a set of general points, with $r\geq 4$, let $\boldsymbol{\mu}=\{\mu_1\geq\ldots\geq\mu_r\}$ be positive integers and let $I:=I(\mathbf{p},\boldsymbol{\mu})$ stand for the corresponding ideal of fat points in $R=k[x,y,z]$.
Then, for any integer  $t\geq \mu_1+\mu_2$ such that $$\frac{(t+2)(t+1)}{2}-\displaystyle\sum_{i=1}^{r}\frac{(\mu_i+1)\mu_{i}}{2}\geq \frac{3\mu_4^2-7\mu_4+4}{2},$$
one has}
$$\dim_k I_t=\frac{(t+2)(t+1)}{2}-\displaystyle\sum_{i=1}^{r}\frac{(\mu_i+1)\mu_{i}}{2}.$$

\medskip

To see how it applies to the above example, note that for this type its fourth multiplicity is $1$, hence the bound prescribed in the above theorem is $0$.
Therefore, we can apply the result in our case with $I:=\widetilde{J}$, for any $t\geq d-2+1=d-1$ such that
$$\frac{(t+2)(t+1)}{2}\geq 
\sum_{i=1}^{r}\frac{(\mu_i+1)\mu_{i}}{2}=\frac{(d-1)(d-2)}{2}+2d-1=\frac{d(d+1)}{2}.
$$
Clearly, any $t\geq d-1$ will do it.
In this range the theorem says that one has an equality
$$\dim_k(\widetilde{J}_t)=\frac{(t+2)(t+1)}{2}-
\frac{d(d+1)}{2}
$$
and  $\dim_k(\widetilde{J}_{d-1})=0$.

This implies that $\dim_k(\widetilde{J}_{t})=0$ for $0<t\leq d-1$ and that the regularity index of $R/\widetilde{J}$ is $d-1$; the latter gives that  $\widetilde{J}$ has no minimal generators in degree $\geq (d-1)+2=d+1$
while in degree $d$ the minimal  number of generators is
$$\dim _k(\widetilde{J}_d)=\frac{(d+1)(d+2)}{2}-\frac{d(d+1)}{2}=d+1.
$$
Therefore, $\widetilde{J}=(\widetilde{J}_d)$ is minimally generated by $d+1$ forms, and hence it must be linearly presented with free resolution of the stated form.

\begin{Example}\rm
Let $T=(d; 8^{r_8},4^{r_4}, 2^{r_2})$ be a homaloidal type, where $d\equiv 1 \pmod{4}$. Then one of the following alternatives takes place:
\begin{enumerate}
\item[$\bullet$] $T=(5;2^6)$ (proper)
\item[$\bullet$] $T=(9;4^4,2^4)$ (proper)
\item[$\bullet$] $T=(13;8,4^6,2^2)$ or $T=(13;8^2,2^{10})$ (both improper)
\item[$\bullet$] $T=(17;8^3,4^6)$ or $T=(17;8^4,2^8)$ (both proper)
\end{enumerate}
\end{Example}
To see this, we first consider the ``halved'' type
\begin{equation}\label{halved_example}
(\left \lceil{d/2}\right \rceil; 4^{r_8},2^{r_4}, 1^{r_2}).
\end{equation}

We could start by crudely bounding $d$ using Noether's inequality; we choose to proceed {\em ab initio} from the equations of condition.

Set $d=4t+1$, so $\left \lceil{d/2}\right \rceil=2t+1$.
Since we are assuming that $T$ is homaloidal, then by Lemma~\ref{doubling} one has the following system of equations for (\ref{halved_example}):
$$\left\{
\begin{array}{ccc}
r_2+2r_4+4r_8 & =& 6t\\
r_2+4r_4+16r_8 &=& 2t(2t+1)
\end{array}
\right.
$$
Subtracting and substituting, yields $r_4, r_2$ as functions of $t$ and $r_8$:

\begin{equation}\label{r_four_and_two}
\left\{
\begin{array}{ccc}
r_4 &=& 2(t(t-1)-3r_8)\\
r_2& = & 2t(5-2t)+8r_8
\end{array}
\right.
\end{equation}
Since the multiplicities are nonnegative, these relations imply the inequalities
$$\frac{t(2t-5)}{4}\leq r_8\leq \frac{t(t-1)}{3},
$$
where the difference between the extremes must be nonnegative.
This forces $t\leq 5$, and hence, $r_8\leq \left \lfloor {20/3}\right \rfloor=6$.
Replacing in the second of the relations in (\ref{r_four_and_two}) now yields $r_2\leq  2t(5-2t)+48$.
Once more, since $r_2\geq 0$, then $t\leq 4$.
From this we have the improved upper bound $r_8\leq 4$.

Finally, the bound $r_8\leq t(t-1)/3$ and the value of $r_2$ also imply that $r_8=0 \Leftrightarrow t\leq 2$.

We now have all the information needed to draw the stated conclusion by a straightforward computation.
The first two on the list correspond, respectively, to the case where $t=1$ and $t=2$ ($r_8=0$).

The case where $t=3$ splits into two sub-cases, according to the two alternatives coming from $r_8\leq t(t-1)/3=2$.

The last item on the list corresponds to $t=4$, where $r_8\leq 4$.
For $r_8=1,2$ the value of $r_2$ would be negative.
Thus, we must have $r_8=3,4$, corresponding to the two laternatives of this item.

\smallskip

To conclude, we check which of the obtained homaloidal types is proper.
The type $(5;2^6)$ is a well-known proper type.
The type $(9;4^4,2^4)$ abuts to $(4;2^3,1^3)$ after composing with two suitable quadratic transformations based on the three highest multiplicities. This is again a well-known proper type.

Next, the first type of degree $13$ ($t=3$) listed above becomes  the type $(7;4^2,\ldots)$ by means of quadratic transformations. The latter is not proper since $7<4+4$.
The second type of degree $13$ listed is not proper since on the face $13<8+8$.

Finally, by quadratic transformations, both types of degree $17$ end up at the type $(4;2^3,1^3)$, hence they are proper.



\noindent {\bf Authors' addresses:}

\medskip

\noindent {\sc Zaqueu Ramos},  Departamento de Matem\'atica, CCET\\ Universidade Federal de Sergipe\\
49100-000 S\~ao Cristov\~ao, Sergipe, Brazil\\
{\em e-mail}: zaqueu@gmail.com\\

\noindent {\sc Aron Simis},  Departamento de Matem\'atica, CCEN\\ Universidade Federal
de Pernambuco\\ 50740-560 Recife, PE, Brazil\\
{\em e-mail}:  aron@dmat.ufpe.br


\begin{thebibliography}{99}

\bibitem{AHT}{I. Aberbach, C. Huneke, N.V. Trung, Reduction numbers, Brian\c{c}on-Skoda theorems and the depth of Rees rings,
Compositio Math. {\bf 97} (1995), 403--434.}

\bibitem{alberich}{M. Alberich-Carrami\~nana, {\sc Geometry of the
Plane Cremona Maps}, Lecture Notes in Mathematics, {\bf 1769},
Springer-Verlag Berlin-Heidelberg,  2002.}


\bibitem{Bocci}{C. Bocci and B. Harbourne, Comparing Powers and Symbolic Powers of Ideals, J. Algebraic Geometry {\bf 19} (2010) 399-–417.}

\bibitem{BHbook}{W. Bruns and J. Herzog, {\sc Cohen--Macaulay Rings}, Cambridge University Press, Cambridge, 1993.}

\bibitem{CiMi}{C. Ciliberto and R. Miranda, Homogeneous interpolation of ten points, J. Algebraic Geometry, {\bf 20} (2011) 685–-726.}


\bibitem{CorZar}{T. Cortadellas and S. Zarzuela, Burch's inequality and the depth of the blow up rings of an ideal,
J. Pure and Applied Algebra {\bf 157} (2001) 183?-204.}


\bibitem{DGM}{E. D. Davis, A. V. Geramita, and P. Maroscia, Perfect homogeneous ideals: Dubreil's theorems revisited, Bull. Sc. Math., 2e s´erie, {\bf 108} (1984), 143--185.}

\bibitem{AHA}{A. Doria, H. Hassanzadeh and A. Simis, A characteristic free criterion of birationality, Advances in Math., {\bf 230} (2012), 390--413.}






\bibitem{Dumnicki5}{M. Dumnicki, T. Szemberg and H. Tutaj-Gasi\'nska, A vanishing theorem and symbolic powers of planar point ideals, LMS J. Comput. Math. {\bf 16} (2013), 373--387.}





\bibitem{FiHaHo}{S. Fitchett, B. Harbourne and S. Holay, Resolutions of Fat Point Ideals involving Eight General
Points of $\pp^2$, J. Algebra, {\bf 244} (2001), 684--705.}

\bibitem{GGP}{A.V. Geramita,A. Gimigliano and Y. Pitteloud, Graded Betti numbers of some embedded rational $n$-folds, Math. Ann., {\bf 301} (1995), 363--380.}

\bibitem{Guardo}{E. Guardo, B. Harbourne and A. Van Tuyl, Symbolic powers versus regular powers of ideals of general points $\pp^1\times \pp^1$,
{\bf 65} (2013), 823--842.}

\bibitem{H0}{B. Harbourne, Complete linear systems on rational surfaces, Trans. Amer, Math. Soc. {\bf 289} (1985), 213–-226.}

\bibitem{H1}{B. Harbourne, The Ideal Generation Problem for Fat Points, J. Pure and Applied Alg. {\bf 145} (2000), 165-–182.}

\bibitem{H2}{B. Harbourne, Problems and Progress: A survey on fat points in P2, Queen’s papers in Pure and Applied Mathematics, The Curves Seminar at Queen’s, vol. 123, 2002.}

\bibitem{HHF}{B. Harbourne, S. Holay and S. Fitchett, Resolutions of ideals of quasiuniform fat point subschemes of $\pp^2,$ Tans. Amer. Math. Soc. 355 (2003), n0. 2, 593--608.}

\bibitem{HS}{S. H. Hassanzadeh and A. Simis, Plane Cremona maps: saturation, regularity and fat ideals, J. Algebra, {\bf 371} (2012), 620--652.}


\bibitem{verde}{J. Hong, A. Simis, W. V. Vasconcelos, Notes on etc., ongoing.}


\bibitem{Nagata1}{M. Nagata, On rational surfaces I, Memoirs of the College of Sciences, University of Kyoto, Series A, {\bf 32} (1960) 351--370.}


\bibitem{Nagata2}{M. Nagata, On rational surfaces II, Memoirs of the College of Sciences, University of Kyoto, Series A, {\bf 33} (1960) 271--293.}

\bibitem{Zaron}{Z. Ramos and A. Simis, Symbolic powers of perfect ideals of codimension 2 and birational maps, J. Algebra 413 (2014) 153--197.}

\bibitem{fat1}{Z. Ramos and A. Simis, Homaloidal nets and ideals of fat points I, LMS Journal of Computation and Mathematics, {\bf 19}  (2016), 54--77.}


\bibitem{ST} {A. Simis and S. Toh\v aneanu, The ubiquity of Sylvester forms in almost complete intersections, Collectanea Math, {\bf 66} (2015) 1--31.}

\bibitem{ram2}{A. Simis,  B. Ulrich and W. V. Vasconcelos, Codimension, multiplicities and integral extensions,
Math. Proc. Camb. Phil. Soc. {\bf 130}, (2001), 237--257.}

\bibitem{Wolmbook3}
W. Vasconcelos, Integral Closure, Rees algebras, multiplicities, algorithms, Springer Monographs in Mathematics, Springer-Verlag, Berlin, 2005.


\end{thebibliography}
\end{document}